\documentclass[11pt]{amsart}
\usepackage{amsbsy}
\usepackage{graphicx,epsfig,subfigure,psfrag}
\usepackage[lmargin=2.1cm,rmargin=2.1cm,tmargin=2.1cm,bmargin=2.9cm]{geometry}

\begin{document}

%%%%%%%%%%%%%%%%%%%%%%%%%%%%%%%%%%%%%%%%%%%%%%%%%%%%%%%%%%%%%%%%%%%%
% Theorem, definition, lemma, proposition, corollary and proof
%%%%%%%%%%%%%%%%%%%%%%%%%%%%%%%%%%%%%%%%%%%%%%%%%%%%%%%%%%%%%%%%%%%%
%%%%%%%%%%%%%%%%%%%%%%%%%%%%%%%%%%%%%%%%%%%%%%%%
\newtheorem{theorem}{Theorem}%[section]
\newtheorem{proposition}{Proposition}%[section]
\newtheorem{lemma}{Lemma}%[section]
\newtheorem{corollary}{Corollary}%[section]%%
\newtheorem{definition}{Definition}%[section]
\newtheorem{remark}{Remark}%[section]
%%%%%%%%%%%%%%%%%%%%%%%%%%%%%%%%%%%
%%%%%%%%%%%%%%%%%%%%%%%%%%%%%%%%%%%%%%%%%%%%%%                  NEW
%%\newcommand{\be}{\begin{equation}}
%%\newcommand{\ee}{\end{equation}}
%%%%%%%%%%%%%%%%%%%%%%%%%%%%%%%%%%%%%%%%%%%%%%%
%%%%%%%%%%%%%%%%%%%%%%%%%%%%%%%%%%%%%%%%%%%%
\newcommand{\tex}{\textstyle}
%%\DeclareMathOperator{\Ind}{Ind}
%% \DeclareMathOperator{\Card}{Card}
%% \DeclareMathOperator{\Deg}{Deg}
%% \DeclareMathOperator{\dist}{dist}
%% \DeclareMathOperator{\Signature}{Signature}
%% \DeclareMathOperator{\MC}{MC}
%% \DeclareMathOperator{\sign}{sign}
%%  \DeclareMathOperator{\Int}{int}
%%%%%%%%%%%%%%%%%%%%%%%%%%%%%%%%%%%%%%%%%%%%
%%%%%%%%%%%%%%%%%%%%%%%%%%%%%%%%%%%%%%%%%%%%
\numberwithin{equation}{section} \numberwithin{theorem}{section}
\numberwithin{proposition}{section} \numberwithin{lemma}{section}
\numberwithin{corollary}{section}
\numberwithin{definition}{section} \numberwithin{remark}{section}
%%%%%%%%%%%%%%%%%%%%%%%%%%%%%%%%%%%%%%%%%%%%
%%%%%%%%%%%%%%%%%%%%%%%%%%%%%%%%%%%%%%%%%%%%
\newcommand{\ren}{\mathbb{R}^N}
\newcommand{\re}{\mathbb{R}}
\newcommand{\n}{\nabla}
\newcommand{\p}{\partial}
\newcommand{\iy}{\infty}
\newcommand{\pa}{\partial}
\newcommand{\fp}{\noindent}
\newcommand{\ms}{\medskip\vskip-.1cm}
\newcommand{\mpb}{\medskip}
%%%%%%%%%%%%%%%%%%%%%%%%%%%%%%%%%%%%%%%%%%%%%%%%%
\newcommand{\AAA}{{\bf A}}
\newcommand{\BB}{{\bf B}}
\newcommand{\CC}{{\bf C}}
\newcommand{\DD}{{\bf D}}
\newcommand{\EE}{{\bf E}}
\newcommand{\FF}{{\bf F}}
\newcommand{\GG}{{\bf G}}
\newcommand{\oo}{{\mathbf \omega}}
\newcommand{\Am}{{\bf A}_{2m}}
\newcommand{\CCC}{{\mathbf  C}}
\newcommand{\II}{{\mathrm{Im}}\,}
\newcommand{\RR}{{\mathrm{Re}}\,}
\newcommand{\eee}{{\mathrm  e}}
%%%%%%%%%%%%%%%%%%%%%%%%%%%%%%%%%%%%%%%%%%%%%%%%%%%%%%%%%%%%%%%%%%%%%%% L^2\rho...
\newcommand{\LL}{L^2_\rho(\ren)}
\newcommand{\LLL}{L^2_{\rho^*}(\ren)}
%%%%%%%%%%%%%%%%%%%%%%%%%%%%%%%%%%
%%%%%%%%%%%%%%%%%%%%%%%%%%%%%%%%%%%%%%%%%%%%%%%%%%%%
\renewcommand{\a}{\alpha}
\renewcommand{\b}{\beta}
\newcommand{\g}{\gamma}
\newcommand{\G}{\Gamma}
\renewcommand{\d}{\delta}
\newcommand{\D}{\Delta}
\newcommand{\e}{\varepsilon}
\newcommand{\var}{\varphi}
\newcommand{\lll}{\l}
\renewcommand{\l}{\lambda}
\renewcommand{\o}{\omega}
\renewcommand{\O}{\Omega}
\newcommand{\s}{\sigma}
\renewcommand{\t}{\tau}
\renewcommand{\th}{\theta}
\newcommand{\z}{\zeta}
\newcommand{\wx}{\widetilde x}
\newcommand{\wt}{\widetilde t}
\newcommand{\noi}{\noindent}
 %%%%%%%%%%%%%%%%%%%%%%%%%%%%%%%%%%%%%%%%%%%
\newcommand{\uu}{{\bf u}}
\newcommand{\xx}{{\bf x}}
\newcommand{\yy}{{\bf y}}
\newcommand{\zz}{{\bf z}}
\newcommand{\aaa}{{\bf a}}
\newcommand{\cc}{{\bf c}}
\newcommand{\jj}{{\bf j}}
\newcommand{\ggg}{{\bf g}}
\newcommand{\UU}{{\bf U}}
\newcommand{\YY}{{\bf Y}}
\newcommand{\HH}{{\bf H}}
\newcommand{\GGG}{{\bf G}}
\newcommand{\VV}{{\bf V}}
\newcommand{\ww}{{\bf w}}
\newcommand{\vv}{{\bf v}}
\newcommand{\hh}{{\bf h}}
\newcommand{\di}{{\rm div}\,}
\newcommand{\ii}{{\rm i}\,}
%%%%%%%%%%%%%%%%%%%%%%%%%%%%%%%%%%
%%%%%%%%%%%%%%%%%%%%%%%%%%%%%%%%%%%%%   VAG, NEW
\newcommand{\inA}{\quad \mbox{in} \quad \ren \times \re_+}
\newcommand{\inB}{\quad \mbox{in} \quad}
\newcommand{\inC}{\quad \mbox{in} \quad \re \times \re_+}
\newcommand{\inD}{\quad \mbox{in} \quad \re}
\newcommand{\forA}{\quad \mbox{for} \quad}
\newcommand{\whereA}{,\quad \mbox{where} \quad}
\newcommand{\asA}{\quad \mbox{as} \quad}
\newcommand{\andA}{\quad \mbox{and} \quad}
\newcommand{\withA}{,\quad \mbox{with} \quad}
\newcommand{\orA}{,\quad \mbox{or} \quad}
\newcommand{\atA}{\quad \mbox{at} \quad}
\newcommand{\onA}{\quad \mbox{on} \quad}
\newcommand{\ef}{\eqref}
\newcommand{\mc}{\mathcal}
\newcommand{\mf}{\mathfrak}

\newcommand{\ssk}{\smallskip}
\newcommand{\LongA}{\quad \Longrightarrow \quad}
%%%%%%%%%%%%%%%%%%%%%%%%%%%%%%%%
%%%%%%%%%%%%%%%%%%%%%%%%%%%%%%%%%%
\def\com#1{\fbox{\parbox{6in}{\texttt{#1}}}}
%%%%%%%%%%%%%%%%%%%%%%%%%%%%%%%%%%
%%%%%%%%%%%%%%%%%%% From Paper1
\def\N{{\mathbb N}}
\def\A{{\cal A}}
\newcommand{\de}{\,d}
\newcommand{\eps}{\varepsilon}
\newcommand{\be}{\begin{equation}}
\newcommand{\ee}{\end{equation}}
\newcommand{\spt}{{\mbox spt}}
\newcommand{\ind}{{\mbox ind}}
\newcommand{\supp}{{\mbox supp}}
\newcommand{\dip}{\displaystyle}
\newcommand{\prt}{\partial}
\renewcommand{\theequation}{\thesection.\arabic{equation}}
\renewcommand{\baselinestretch}{1.1}
%%%%%%%%%%%%%%%%%%%%%%%%%%%%%%%%%%%%%%%%%%%%%%%
\newcommand{\Dm}{(-\D)^m}

%%%%%%%%%%%%%%%%%%%%%%%%% VICTOR
\title
%%%%%[Positivity]
%%%%%%%%%%%%%%%%%%%%%%%%%
%%%%{\bf How  $\sqrt{\mbox{log\,log}}$ factor occurs in blow-up
{\bf The Cauchy problem for  tenth-order thin film equation I.
Bifurcation of oscillatory fundamental solutions}

\author{P.~\'Alvarez-Caudevilla, J.D. Evans and V.A.~Galaktionov}

\address{Universidad Carlos III de Madrid, Spain}
\email{pacaudev@math.uc3m.es}

\address{Department of Mathematical Sciences, University of Bath,
 Bath BA2 7AY, UK}
%%% and Keldysh Institute of Applied Mathematics,
%%%% Miusskaya Sq. 4, 125047 Moscow, RUSSIA}
\email{masjde@bath.ac.uk}

\address{Department of Mathematical Sciences, University of Bath,
 Bath BA2 7AY, UK}
%%% and Keldysh Institute of Applied Mathematics,
%%%% Miusskaya Sq. 4, 125047 Moscow, RUSSIA}
\email{vag@maths.bath.ac.uk}

 

\thanks{This works has been partially supported by the Ministry of Economy and Competitiveness of
Spain under research project MTM2012-33258.}
\date{\today}

%%%%%%%%%%%%%%%%%%%%%%%%%%%%%%%%%%%%%%%%%%%%%%%%%%%%%%%%

%%\setlength{\topmargin}{5mm} \pagestyle{myheadings} \markboth{}
%%{Tenth-order thin film equation}

%%\begin{document}

%%%%$$ $$

\begin{abstract}

 Fundamental global similarity solutions of the tenth-order thin
 film equation
 $$ %%%\begin{equation}
 %%%\label{i1a}
    u_{t} = \nabla \cdot(|u|^{n} \n \D^4 u)
 \quad \mbox{in} \quad \ren \times \re_+
    \,,
  $$ %%%%  \ee
    where $n>0$ are studied. The main approach consists in passing
    to the limit $n \to 0^+$ by using Hermitian non-self-adjoint
    spectral theory corresponding to the rescaled linear
    poly-harmonic equation
     $$u_t= \D^5 u\quad \mbox{in} \quad \ren \times \re_+
    \,.
  $$
%%%\end{equation}

\end{abstract}

%%%%%%%%%%%%%%%%%%%%%%%%%%%
\maketitle

%%%%%%%%%%%%%%%%%%%%%%%%%%%

%%%%%%%%%%%%%%%%%%%%%%%%%%%%%%%%%%%%%%%%%%%%%%%%
\section{Introduction: the TFE-10 and  nonlinear eigenvalue problem }
 \label{S1}

\subsection{Main model and result: toward discrete real nonlinear spectrum}

\noindent
    We study the global-in-time behaviour of  solutions of
the tenth-order quasilinear evolution equation of parabolic type,
called the {\em thin film equation} (TFE--10)
\begin{equation}
\label{i1}
    u_{t} = \nabla \cdot(|u|^{n} \n \D^4 u)
 \quad \mbox{in} \quad \ren \times \re_+
    \,,
\end{equation}
%%\com{PAC: He cambiado el operador $\n^9$ aqui y despues, pero vas
%%%a verificar...}
where $\n={\rm grad}_x$ and $n>0$ is a real parameter. In view of
the degenerate mobility coefficient $|u|^n$, equation \eqref{i1}
is written for solutions of changing sign, which can occur in the
{\em Cauchy problem} (CP) and also in some {\em free boundary
problems} (FBPs).

Equation \ef{i1} has been chosen as a typical higher-order
quasilinear degenerate parabolic model, which is very difficult to
study, and this is key for us; see below. Although the fourth-order version
has been the most studied, the sixth-order TFE is known to occur
in several applications and, during the last ten-fifteen years,
has begun to steadily penetrate into modern nonlinear PDE theory; see
references in \cite[\S~1.1]{EGK3} and more recently \cite{CG,Liu1,Liu,LQ} and 
\cite{GuiLon,GB} where several applications of these problems are described, in particular 
image processing. 

However, our main intention here is to develop the mathematical theory in the analysis
of degenerate even higher-order equations without looking at the applications which we are not aware of up to order eighth. The analysis
performed in this work will provide some new techniques in obtaining qualitative results for these difficult to analyze PDEs. 
Since there has been a lot of published material about fourth and sixth order higher order equations of similar form to \eqref{i1} 
we have jumped to tenth order to generalize this theory for even higher-order degenerate equations of this type.
%%%% JDE to Pablo -
% I've removed reference to 8th order - What is the application in which it arises ?
% I'm not aware of any - if you can find a reference for an application - then stick it back in with ref.
%
% I've put in some more recent refs on 6th order - Can you have a look as well for any new references/work on
% 6th order thin-film/porous medium equation
% I refereed a paper by Liu and Tian Weak solutions for a sixth-order porous medium type equation
% can't find it though - may be not published ? (I was positive on it !)
%%%%%%%%%%%%%

Let us state our main result. In Section \ref{S2}, we introduce
global self-similar solutions of \ef{i1} of the standard form
\begin{equation}
\label{sf3}
 \tex{
    u(x,t):= t^{-\a} f(y), \quad \hbox{with}\quad
    y:=\frac{x}{t^\b},\quad \b= \frac{1-n \a}{10},
    }
\end{equation}
 where $f$ satisfies an elliptic equation given below.
 Then a  {\em nonlinear eigenvalue problem} with a nonlinear
 {\em real eigenvalue} $\a$ occurs\footnote{More precisely, since the eigenvalue
 $\a$ enters not only the standard term $\a f$, but also the linear differential one
 $\frac{1-\a n}{10}\, y \cdot \nabla f$, it is more correct to talk about a
 ``linear (in $\a$) spectral pencil for the quasilinear TFE-10 operator",
 though, for simplicity, we keep referring to the nonlinear eigenvalue problem.
 In contract to these nonlinear issues, for $n=0$,
the second term looses $\a$, and we arrive a standard linear
eigenvalue problem for the non-self-adjoint operator ${\bf B}=\D^5
+ \frac 1{10}\, y \cdot \n + \frac N{10} \, I$; see Section
\ref{S3}.}:
\begin{equation}
\label{self1}
 %%\fbox{$
 \tex{
     \nabla \cdot \left[ |f|^{n} \n \D^4 f\right] +\frac{1-\a n}{10}\, y \cdot \nabla f +\a
    f=0,
    \quad f \in C_0(\ren)\, ,
    }
   %% $}
\end{equation}
 where the problem setting includes finite propagation phenomena
 for such TFEs, i.e., $f$ is assumed to be compactly supported, $f \in
 C_0(\ren)$. This is a kind of an assumed ``minimal"
behaviour of $f(y)$ as $y \to \iy$.

Using long-established terminology, we call such similarity
solutions \ef{sf3} (and also the corresponding profiles $f_\g$) to
be a sequence of  {\em fundamental solutions}. Though, actually,
the classic fundamental solution is the first radially symmetric
one (with the first kernel $f_0=f_0(|y|)$), which is  the {\em instantaneous
source-type solution} of \ef{i1} with Dirac's delta as initial data. Moreover, for
$n=0$, $f_0(|y|)$ becomes the actual rescaled kernel of the fundamental solutions  of the linear operator
$D_t-\D_x^5$.

\ssk

Our main goal is to show {\em analytically} that, at least, for small $n>0$,
 \be
 \label{main1}
 %% \fbox{
 \mbox{(\ref{self1}) admits a countable set of fundamental
 solutions $\Phi(n)=\{\a_\g,f_\g\}_{|\g| \ge 0}$},
  %%}
  \ee
  where $\g$ is a multiindex in $\ren$ to numerate   these eigenvalue-eigenfunction pairs. Global extensions
  of such ``$n$-branches" of some first fundamental solutions can be checked numerically.

%%%%%%%%%%%%%%%%%%%%%%%%%%%%%%%%%%%%%%%%%%%%%%%%%%%
\subsection{First discussion: possible origins of discrete
nonlinear spectra and principle difficulties}

It is key for us that \ef{self1} {\em is not variational}, then
we cannot use powerful tools such as Lusternik--Schnirel'man (L--S, for
short) category-genus theory, which in many cases is known to
provide a {\em countable} family of critical points (solutions), if
the category of the functional subset involved is infinite.

It is crucial, and well known, that the L--S min-max approach
{\em does not detect all families of critical points}. However,
sometimes it can revive a minor amount of solutions. A somehow
special example was revealed   in \cite{GMPSobI, GMPSobII}, where
key features of those variational L--S and fibering approaches are
available. Namely,  for some variational fourth-order and
higher-order ODEs in $\re$, including  those with the typical
nonlinearity $|f|^n f$, as above,
  \be
 \label{mm.561}
  \mbox{$
 -  (|f|^n f)^{(4)} + |f|^n f= \frac 1n \, f \inB  \re, \quad f \in C_0(\re) \quad (n>0),
 $}
  \ee
 as well as for the following standard looking
one with the only cubic nonlinearity \cite[\S~6]{GMPSobII}:
 \be
 \label{anal1}
-f^{(4)}+ f = f^3 \quad \mbox{in} \quad \re, \quad f \in
H^4_\rho(\re) \quad (\rho={\mathrm e}^{a|y|^{4/3}}, \,\, a>0
\,\,\, \mbox{small}),
 \ee
  it was shown that these admit a {\em countable set of
countable families of solutions}, while the L--S/fibering approach
detects only {\sc one} such discrete family of (min-max) critical
points. Further countable families are not expected to be
determined easily by more advanced techniques of potential theory,
such as the mountain pass lemma, and others. Existence of other,
not L--S type critical points for \ef{mm.561} and \ef{anal1} were
shown in \cite{GMPSobI, GMPSobII} by using a combination of
numerical and (often, formal) analytic methods and heavy use of
oscillatory nature of solutions close to finite interfaces (for
\ef{mm.561}) and at infinity (for \ef{anal1}). In particular,
detecting the corresponding L--S countable sequence of critical
points was done {\em numerically}, i.e., by checking their actual
min-max features (their critical values must be maximal among
other solutions belonging to the functional subset of a given
category, and having  a ``suitable geometric shape").

 Therefore, even in
the variational setting, counting various families of critical
points and values represents a difficult open problem for such
higher-order ODEs, to say nothing of their elliptic extensions in
$\ren$.

Hence, studying the nonlinear eigenvalue problem \ef{self1}, we
will rely on a different approach, which is also effective for
such difficult variational problems and detects more solutions
than L--S/fibering theory (though only locally upon the
parameter). Namely,
 our main approach is the idea of a
``homotopic deformation" of \ef{i1}  as $n \to 0^+$ (Section
\ref{S4}) and reducing it to the classic {\em poly-harmonic
equation of tenth order}
\begin{equation}
\label{lin5}
    \tex{ u_{t} = \D^5 u\quad \hbox{in} \quad \re^{N} \times \re_+\,.}
\end{equation}
The corresponding \ef{self1} then reduces to a standard (but not
self-adjoint) Hermitian-type  linear eigenvalue problem, which is
 treated in Section \ref{S3}.
 Therefore, according to this approach, the nonlinear version of
 \ef{main1} has the origin in the discreteness-reality of the
 spectrum of the corresponding linear operator.

 Finally, in Section \ref{S5}, we present numerical results for eigenfunctions with explicit eigenvalues. These are 
 the eigenfunctions in the $n=0$ case, which provide the starting points of the n-branch solutions. The eigenfunctions in the mass 
 conservative case are also presented for selected $n$, which constitutes the first n-branch. 
%These provide
% confirming that
% first bifurcation $n$-branches of solutions (locally studied above close to $n=0$) admit a clear global
% extentions in $n$, up to $n=100$ and more.
% \com{JDE: Please correct, this is just to recall about this new Section. Put extra details, anything.}

%%%%%%%%%%%%%%%%%%%%%%%%%%%%%%%%%%%
\subsection{The second model: bifurcations in $\re^2$}

 In Appendix A, we show how to extend our homotopy approach to a
 more complicated   {\em unstable thin film equation} (TFE--10) in
 the critical case
\begin{equation}
\label{e1}
 \tex{
    u_{t} = \nabla \cdot(|u|^{n} \n \D^4 u)-\D(|u|^{p-1}u)
 \quad \mbox{in} \quad \ren \times \re_+
    \,, \quad p>n+1,%%%%%\quad p=n+1+ \frac 8N,
    }
\end{equation}
with the extra unstable diffusion term. We briefly and formally
show that, revealing a discrete real nonlinear spectrum for
\ef{e1} then requires a simultaneous {\em double} homotopy
deformation $n \to 0^+$ and $p \to 1^+$ leading to a new linear
Hermitian spectral theory. We do not develop it here and just
focus on a principal opportunity to detect a discrete nonlinear
spectrum for \ef{e1}.

%%%\ssk

%%%%%%%%%%%%%%%%%%%%%%%%%%%%%%%%%%%%%%%%%
\subsection{Global extension of bifurcation branches: a principal
open problem}

It is worth mentioning that, for both problems \ef{self1} and the
corresponding problem occurring for \ef{e1} (after the similarity
time-scaling; see \ef{sfe5}), global extension of bifurcation
$n$-branches ($(n,p)$-branches for \ef{e1}) represents a difficult
open problem of general nonlinear operator theory. Moreover, as
was shown in \cite{GalPetII} (see also other examples in
\cite{GMPSobII}), the TFE-4 with absorption $-|u|^{p-1}u$ (instead
of the backward-in-time diffusion as in \ef{self1}), depending on
not that small $n \sim 1$, has some $p$-bifurcation branches can
have turning (saddle-node) points and can represent a closed
loops, so that such branches are not globally extended. On the
other hand, for equations with monotone operators such as the
PME-4 (see \ef{PME4} below), the $n$-branches seem to be globally
extensible in $n>0$, \cite{GalRDE4n}.

%%%%%%%%%%%%%%%%%%%%%%%%%%%%%%%%%%%%%%
\subsection{Back to our main motivation}

After posing our main models to study, we must confess that our
main motivation to chose those was their actual extreme
mathematical difficulty. We wanted to see which mathematical
methods and ideas could be applied to justify \ef{main1} using
{\em any kind} of mathematical tools.

 Though we were not able to
justify our results rigorously (and we suspect that this cannot be
done in principle), we believe that our homotopic approach is the
only one available for declaring the result \ef{main1}, which, as
we claim, is in fact a generic property of many nonlinear
eigenvalue problems for elliptic equations. Indeed, we also claim
that the {\em discreteness}  of the nonlinear spectrum in
\ef{main1} and the {\em reality} of all the eigenvalues have their
deep roots in the linear Hermitian spectral theory corresponding
to $n=0$. Thus, we observe how a {\em ``nonlinear spectral theory
bifurcates from a linear non-self-adjoint one"}.
%%\footnote{If
%%the Readers allow us to say like that.}.

\ssk

 Note that the elliptic
equation \ef{self1} with an extra parameter $\a$ is very
difficult to analyse even in one-dimension, where it becomes a tenth-order
ODE nonlinear eigenvalue problem:
 \be
 \label{ODE1}
  \tex{
 (|f|^n f^{(9)})' +\frac{1-\a n}{10}\, y \, f' +\a
    f=0,
    \quad f \in C_0(\re)\, .
    }
    \ee
 Indeed, this ODE creates  a 10-dimensional phase space and a construction of
 suitable homotopic connections of equilibria, which are
 admissible for necessary nonlinear eigenfunctions, is not easy at
 all. In the forthcoming paper \cite{AEG2}, we study the first
 eigenfunction of \ef{ODE1}, i.e., the fundamental source-type
 profile $f_0(y)$ by using a variety of other analytical,
 asymptotic, and numerical methods.

%%%%%%%%%%%%%%%%%%%%%%%%%%%%%%%%%%%%%%%%%
\setcounter{equation}{0}
\section{Problem setting and self-similar solutions}
\label{S2}

%%%%%%%%%%%%%%%%%%%%%%%%%%%%%%%%%%%

%%%%%%%%%%%%%%%%%%%%%%%%%%%%%%%%%%%
\subsection{The FBP and CP}

As earlier in \cite{EGK1}--\cite{EGK4}, we distinct the standard
{free-boundary problem} (FBP) for \ef{i1} and the {\em Cauchy
problem}; see further details therein.

 For  both the FBP and the
CP, the solutions are assumed to satisfy standard free-boundary
conditions or boundary conditions at infinity:
\begin{equation}
\label{i3}
    \left\{\begin{array}{ll} u=0, & \hbox{zero-height,}\\
    \nabla u=\nabla^2 u=\nabla^3 u=\nabla^4 u=0,\\
    -{\bf n} \cdot (|u|^{n} \nabla \Delta^4 u)=0, &
    \hbox{conservation of mass (zero-flux)}\end{array} \right.
\end{equation}
at the singularity surface (interface) $\Gamma_0[u]$, which is the
lateral boundary of
\begin{equation}
 \label{gamma1}
    \hbox{supp} \;u \subset \ren \times \re_+,\quad N \geq 1\,,
\end{equation}
where ${\bf n}$ stands for the unit outward normal to
$\Gamma_0[u]$.  Note that, for sufficiently smooth interfaces,
%%%, for the FBP
the condition on the flux can be read as
\begin{equation*}
    \lim_{\hbox{dist}(x,\Gamma_0[u])\downarrow 0}
    -{\bf n} \cdot \nabla (|u|^{n}  \Delta^4 u)=0.
\end{equation*}

\ssk

For the CP, the assumption of nonnegativity is got rid of, and
solutions become oscillatory close to interfaces. It is then key,
for the CP, that the solutions are expected to be
``smoother" at the interface than those for the FBP, i.e., \ef{i3}
are not sufficient to define their regularity. These {\em maximal
regularity} issues for the CP, leading to oscillatory solutions,
are under scrutiny in \cite{EGK2} for a fourth-order case. However, since as far as we know there is no
knowledge of how the solutions for these problems should be, little more can be said about it.

Moreover, we denote by
\begin{equation}
 \label{Mass1}
 \tex{
    M(t):=\int u(x,t) \, {\mathrm d}x
    }
\end{equation}
the mass of the solution, where integration is performed over
smooth support ($\ren$ is allowed for the CP only). Then,
differentiating $M(t)$ with respect to $t$ and applying the
divergence theorem (under natural regularity assumptions on
solutions and free boundary), we have that
\begin{equation*}
 \tex{
  J(t):=  \frac{{\mathrm d}M}{{\mathrm d}t}= -
  \int\limits_{\Gamma_0\cap\{t\}}{\bf n} \cdot
     (|u|^{n} \nabla \Delta^4 u)\, .
     }
\end{equation*}
The mass is conserved if
   $ J(t) \equiv 0$, which is assured by the flux condition in
   \eqref{i3}.
The problem is completed with bounded, smooth, integrable,
compactly supported initial data
\begin{equation}
\label{i4}
    u(x,0)=u_0(x) \quad \hbox{in} \quad \Gamma_0[u] \cap \{t=0\}.
\end{equation}
In the CP for \eqref{i1} in $\ren \times \re_+$, one needs to pose
bounded compactly supported initial data \eqref{i4} prescribed in
$\ren$. Then,  under the same zero flux condition at finite
interfaces (to be established separately), the mass is preserved.

%%%%%%%%%%%%%%%%%%%%%%%%%%%%%%%%%%%%%%%%%%%%%%%%%%%%%
\subsection{Global similarity solutions: towards a nonlinear eigenvalue problem}

 We now begin to specify the self-similar
solutions of the equation \eqref{i1}, which
are admitted due to its natural scaling-invariant  nature.
In the case of the mass being conserved, we have global in time
source-type solutions.

The equation \eqref{i1} is invariant under the two-parameter scaling group
\begin{equation*}
    x:= \mu \bar x,\quad t:= \l \bar t,\quad u:= \left(\frac{\mu^{10}}{\lambda} \right)^{\frac{1}{n}} \bar u.
\end{equation*}
Taking a power law dependence $\mu=\lambda^\beta$, motivates the consideration of self-similar solutions in the form
\begin{equation*}
\tex{
    u(x,t):= \l^{\frac{10\b-1}{n}} f(\frac{x}{t^\b}),
    }
\end{equation*}
as in \ef{sf3}.
%%with
%\begin{equation*}
% \tex{
%    \frac{\p u}{\p t}= \frac{\nu}{\l} \frac{\p \bar u}{\p \bar t},\quad
%    \frac{\p u}{\p x_{i}}= \frac{\nu}{\mu} \frac{\p \bar u}{\p \bar{x}_{i}},\quad
%    \frac{\p^k u}{\p x_{i}^k}= \frac{\nu}{\mu^k} \frac{\p^k \bar u}{\p
%    \bar{x}_{i}^k}, \quad \hbox{with}\quad k=1,\cdots,9.
%     }
%\end{equation*}
% and substituting those expressions in \eqref{i1} yields
%\begin{equation*}
% \tex{
%    \frac{\nu}{\l} \frac{\p \bar u}{\p t}=
%    \frac{\nu^{n+1}}{\mu^{10}} \nabla \cdot (|\bar u|^{n} \n \D^4 \bar u)\,.
%    }
%\end{equation*}
%To keep this equation invariant, the following  must be fulfilled:
%\begin{equation}
%\label{sf2}
% \tex{
%    \frac{\nu}{\l}=\frac{\nu^{n+1}}{\mu^{10}}
%     %%% \quad \mbox{so that}
%   %%% }
%%%%\end{equation}
%%%\begin{equation*}
%  \,\,\Rightarrow \,\,  \mu := \l^\b \Longrightarrow \nu := \l^{\frac{10\b-1}{n}}, \quad
%   %%% \mbox{and}\quad
% %%%%\tex{
%    u(x,t):= \l^{\frac{10\b-1}{n}} \bar u(\bar x,\bar t) = \l^{\frac{10\b-1}{n}}
%    \bar u(\frac{x}{\mu},\frac{t}{\l}).
%    }
%\end{equation}
%Consequently,
%\begin{equation*}
% \tex{
%    u(x,t):= t^{\frac{10\b-1}{n}} f(\frac{x}{t^\b}),
%     }
%\end{equation*}
%where $t=\l$ and $f(\frac{x}{t^\b})=\bar u(\frac{x}{t^\b},1)$.
%Consequently, comparing to \eqref{sf2}, we obtain
%\begin{equation*}
%   n\a +10\b=1,
%\end{equation*}
%which links the parameters $\a$ and $\b$ as in \ef{sf3}.
Hence,
substituting \ef{sf3} into \eqref{i1} and rearranging terms, we
find that the function $f$ solves the quasilinear elliptic equation
given in \ef{self1}.
%% of the form
%%\begin{equation}
%%\label{sf4}
%%    \nabla \cdot \left[ |f|^{n} \n \D^4 f\right] +\b \,y \cdot \nabla f+\a f=0\,.
%%\end{equation}
We thus finally arrive at the {nonlinear eigenvalue problem} \ef{self1},
%%\begin{equation}
%%\label{sf5}
%% \fbox{$
%% \tex{
%%     \nabla \cdot \left[ |f|^{n} \n \D^4 f\right] +\frac{1-\a n}{10}\, y \cdot \nabla f +\a
%%    f=0,
%%    \quad f \in C_0(\ren)\, ,
%%    }
%%    $}
%%\end{equation}
where we add to the elliptic equation  a natural assumption that
 $f$ must be compactly supported (and, of course, sufficiently
 smooth at the interface, which is an accompanying question to be
 discussed as well).

Thus, for such degenerate elliptic equations,
  the functional setting of \ef{self1} assumes that we are
 looking for  (weak) {\em compactly supported} solutions $f(y)$ as
 certain ``nonlinear eigenfunctions" that hopefully occur for special values of nonlinear eigenvalues
  $\{\a_\g\}_{|\g| \ge 0}$. Therefore,  our goal is to justify
  that
 \ef{main1} holds.

  %%, labelling  the eigenfunctions via a multiindex $\s$,
 %%\be
 %%\label{sf51}
 %%\fbox{$
 %%\mbox{
 %%(\ref{self1}) possesses  a countable set of
 %%eigenfunction/value pairs $\{f_k,\, \a_k\}_{|\s|=k \ge 0}$.
 %% }
 %% $}
 %% \ee

 Concerning  the well-known properties of finite propagation for TFEs, we refer to papers
 \cite{EGK1}--\cite{EGK4}, where a large amount of earlier
 references are available; see also \cite{GMPSobI, GMPSobII} for more recent
 results and references in this elliptic area.
 However, one should observe that there are still
a few entirely rigorous results, especially those that are
attributed to the Cauchy problem for TFEs.

In the linear case $n=0$,
 the condition $f \in C_0(\ren)$, is naturally replaced by the requirement that the
 eigenfunctions $\psi_\b(y)$ exhibit typical exponential decay at
 infinity, a property that is reinforced by introducing  appropriate weighted $L^2$-spaces. Complete details about the spectral theory
 for this linear problem when $n=0$ shortly.
Actually,
 using the homotopy limit $n \to 0^+$, we will be obliged for
 small $n>0$,
 instead of $C_0$-setting in
\eqref{self1}, to use the following weighted $L^2$-space:
  \begin{equation}
   \label{WW11}
  f \in L^2_\rho(\ren), \quad \mbox{where} \quad \rho(y)={\mathrm e}^{a |y|^{10/9}},
  \quad a>0 \,\,\,\mbox{small}.
  \end{equation}

Note that, in the case of the Cauchy problem with conservation of
mass making use of the self-similar solutions \eqref{sf3}, we have
that
\begin{equation}
 \label{Mass11}
 \tex{
    M(t):=\int\limits_{\ren} u(x,t) \, {\mathrm d}x=t^{-\a} \int\limits_{\ren}
     f\big(\frac{x}{t^\b}\big)
    \, {\mathrm d}x = t^{-\a+\b N} \int\limits_{\ren} f(y)
    \, {\mathrm d}y,
    }
\end{equation}
 where the actual integration is performed over the support
  ${\rm supp}\, f$
of the nonlinear eigenfunction.
 Then, as is well known,
if $\int f \not = 0$,
  the exponents are calculated giving the first explicit nonlinear eigenvalue:
\begin{equation}
\label{alb1}
 \tex{
 -\a + \b N=0 \LongA
    \a_0(n)=\frac{N}{10+Nn} \quad \mbox{and}  \quad \b_0(n)=\frac{1}{10+Nn}.
    }
\end{equation}

%%%%%%%%%%%%%%%%%%%%%%%%%%%%%%%%%%%%%%%%%%%%%%%%%%%%%
\setcounter{equation}{0}
\section{Hermitian spectral theory of the linear rescaled operators}
\label{S3}

The Hermitian spectral theory developed in
\cite{EGKP} for a pair $\{\BB,\BB^*\}$ of linear rescaled
operators  for $n=0$, i.e., for the {\em poly-harmonic equation}
\begin{equation}
\label{s1}
 \tex{
    u_{t} = \D^{5} u\quad \hbox{in} \quad \ren \times \re_+\,,
    \quad \BB=\D^{5} +\frac{1}{10}\, y \cdot \n+ \frac{N}{10}\, I, \quad
    \BB^*=\D^{5} - \frac{1}{10} \, y \cdot \n,
    }
\end{equation}
whose solutions are $C^\infty$, have infinite speed of propagation
and oscillate infinitely near the interfaces will be essential for our further analysis to consider the (homotopic) limit $n \to 0^+$ for having
a better understanding of the singular
oscillatory properties of the solutions of the CP for (\ref{i1}). Therefore, in this section,
we establish the spectrum $\s(\bf{B})$
of the linear operator $\bf{B}$ obtained from the rescaling of the
linear counterpart of \eqref{i1}, i.e., the poly-harmonic equation of tenth order.

%%%%%%%%%%%%%%%%%%%%%%%%%%%%%%%%%%%%%
\subsection{How the operator $\BB$ appears: a linear eigenvalue problem}

Let $u(x,t)$ be the unique solution of the CP for the linear
parabolic poly-harmonic equation of tenth order \eqref{s1} with
the initial data (the space as in \eqref{WW11} to be more properly
introduced shortly)
 $ %%\begin{equation*}
    u_0 \in L_{\rho}^2(\ren),
$ %%%\end{equation*}
given by the convolution Poisson-type integral
\begin{equation}
\label{s2}
 \tex{
    u(x,t)=b(t)\, * \, u_0 \equiv t^{-\frac{N}{10}} \int\limits_{\ren} F((x-z)t^{-\frac{1}{10}})
     u_0(z)\, {\mathrm d}z.
    }
\end{equation}
Here, by scaling invariance of the problem, in a similar way as
was done in the previous section for \eqref{i1}, the unique
fundamental solution of the operator $\frac{\p}{\p t} - \D^{5}$ has
the self-similar structure
\begin{equation}
\label{s3}
 \tex{
    b(x,t)=t^{-\frac{N}{10}} F(y), \quad y:=\frac{x}{t^{1/10}} \quad  (x\in \ren).
    }
\end{equation}
Substituting $b(x,t)$ into \eqref{s1} yields that the rescaled
fundamental kernel $F$ in \ef{s3} solves the linear elliptic
problem
\begin{equation}
\label{s4}
 \tex{
    {\bf B}F \equiv \D_y^{5} F + \frac{1}{10}\, y \cdot \nabla_y F  +\frac{N}{10}\, F=0
    \quad \hbox{in} \quad \ren,\quad \int\limits_{\ren} F(y) \, {\mathrm
    d}y=1.
    }
\end{equation}

${\bf B}$ is a non-symmetric linear operator, which is bounded
from $H_{\rho}^{10}(\ren)$ to $L_{\rho}^2(\ren)$ with the exponential
weight as in \ef{WW11}.
Moreover,  $a\in (0,2d)$ is any positive constant, depending on the
parameter $d>0$, which characterises  the exponential decay of the
kernel $F(y)$:
 \begin{equation}
  \label{F11}
 |F(y)| \le D \, {\mathrm e}^{-d |y|^{10/9}} \quad \mbox{in} \quad
 \ren,%%%% \quad \big(D>0, \,\,\, d=3 \cdot 2^{-\frac{11}3}\big).
   \end{equation}
   where $D>0$ is a constant and $d$ is the maximal negative real
   part of roots of the  equation
    $$
    \tex{
    a^9=- \frac 1{10} \big( \frac 9{10}\big)^9.
    }
    $$

 %%\com{VAG: is that d correct?}

%% \com{PAC: no, corrected, have a look...}

 By  $F$ we denote the oscillatory rescaled kernel as
the only solution of \eqref{s4}, which has exponential decay,
oscillates as $|y|\rightarrow \infty$, and satisfies the  standard
 estimate \ef{F11}.

\ssk

Thus, we need to solve the corresponding {\em linear eigenvalue
problem}:
 \be
 \label{LP1}
  \fbox{$
  \BB \psi = \l \psi \inB \ren, \quad \psi \in H^{10}_\rho(\ren).
   $}
   \ee
   One can see that the nonlinear
   one (\ref{self1}) formally reduces to (\ref{LP1}) at $n=0$ with
   the following shifting of the corresponding eigenvalues:
    $$
    \tex{
    \l=-\a + \frac{N}{10}.
    }
    $$
 In fact, this is the main reason to calling (\ref{self1}) a
 nonlinear eigenvalue problem, and, crucially, the discreteness of
 the real spectrum of the linear one (\ref{LP1}) will be shown to
 be inherited by the nonlinear problem, but we are still a long way from
 justifying such an issue.

%%%%%%%%%%%%%%%%%%%%%%%%%%%%%%%%%%%%%%
\subsection{Functional setting and semigroup expansion}

    Thus,
 we solve (\ref{LP1}) and calculate the spectrum of $\s({\bf B})$ in the
weighted space $L_{\rho}^2(\ren)$. We then need the following
Hilbert space:
%%% such that
\begin{equation*}
    H_{\rho}^{10}(\ren) \subset L_{\rho}^2(\ren) \subset L^2(\ren).
\end{equation*}
The Hilbert space $H_{\rho}^{10}(\ren)$ has the following inner product:
\begin{equation*}
 \tex{
    \big\langle v,w \big\rangle_{\rho} := \int\limits_{\ren} \rho(y) \sum\limits_{k=0}^{10}
     D^{k} v(y) \overline{{D^{k} w(y)}}
     \, {\mathrm d}y,
      }
\end{equation*}
where $D^{k} v$ stands for the vector
    $\{D^{\b} v\,,\,|\b|=k\}$,
and the norm
\begin{equation*}
 \tex{
    \| v \|_{\rho}^2 := \int\limits_{\ren} \rho(y) \sum\limits_{k=0}^{10}
     |D^{k} v(y)|^2 \, {\mathrm d}y.
    }
\end{equation*}

\par
 Next, introducing the rescaled variables
\begin{equation}
\label{s6}
 \tex{
    u(x,t)=t^{-\frac{N}{10}} w(y,\tau), \quad y:=\frac{x}{t^{1/10}}, \quad \tau= \ln t \,:\,
    \re_+ \to \re,
  }
\end{equation}
we find that the rescaled solution $w(y,\t)$ satisfies the
evolution equation
\begin{equation}
\label{s7}
    w_{\tau} = {\bf B}w\,,
\end{equation}
since, substituting the representation of $u(x,t)$ \eqref{s6} into
\eqref{s1} yields
\begin{equation*}
 \tex{
    \D_y^{5} w + \frac{1}{10} \, y \cdot \nabla_y w  +\frac{N}{10} \, w= t\, \frac{\p w}{\p t} \frac{\p \tau}{\p t}.
  }
\end{equation*}
Thus, to keep this invariant, the following should be satisfied:
\begin{equation*}
 \tex{
    t \, \frac{\p \tau}{\p t}=1 \,\, \Longrightarrow \,\, \tau = \ln t, \quad \mbox{i.e.,
    as defined in (\ref{s6})}.
    }
\end{equation*}
Hence, $w(y,\tau)$ is the solution of the Cauchy problem for the
equation \eqref{s7} and with the following initial condition at
$\tau=0$, i.e., at $t=1$:
\begin{equation}
\label{s8}
    w_0(y) = u(y,1)\equiv b(1)\, * \, u_0 = F\, * \, u_0\, .
\end{equation}
Then, the linear operator $\frac{\p}{\p \tau} - {\bf B}$ is also a
rescaled version of the standard parabolic one $\frac{\p}{\p t} -
\D^{5}$. Therefore, the corresponding semigroup ${\mathrm e}^{{\bf
B} \tau}$ admits an explicit integral representation. This helps
to establish some properties of the operator ${\bf B}$ and
describes other evolution features of the linear flow.

Indeed, from \eqref{s2} we find the following explicit representation of the
semigroup:
\begin{equation*}
%%\label{s9}
 \tex{
    w(y,\tau)=\int\limits_{\ren} F \big(y-z{\mathrm e}^{-\frac{\tau}{10}}\big)\, u_0(z) \,
    {\mathrm d}z \equiv {\mathrm e}^{{\bf B} \tau} w_0, \quad \mbox{where}
     \quad
    x=t^{\frac{1}{10}}y,  \quad \tau=\ln t.
    }
\end{equation*}
Subsequently, consider Taylor's power series of the analytic
kernel\footnote{We hope that returning here to the multiindex $\b$
instead of $\s$ in \ef{self1} will not lead to a confusion with
the exponent $\b$ in self-similar scaling \ef{sf3}.}
\begin{equation}
\label{s10}
 \tex{
    F\big(y-z {\mathrm e}^{-\frac{\tau}{10}}\big)=\sum\limits_{(\b)} {\mathrm e}^{
    -\frac{|\b|\tau}{10}} \frac{(-1)^{|\b|}}{\b!}    D^\b F(y) z^\b
    \equiv \sum\limits_{(\b)} {\mathrm e}^{-\frac{|\b|\tau}{10}} \frac{1}{\sqrt{\b!}} \psi_\b(y) z^\b,
    }
\end{equation}
for any $y\in \ren$, where
%%\begin{equation*}
    $z^{\b}:=z_1^{\b_1}\cdots z_{N}^{\b_{N}}$
%%\end{equation*}
and $\psi_{\b}$ are the normalized eigenfunctions for the operator
$\bf{B}$.
The series in \ef{s10} converges uniformly on compact subsets in
$z \in \ren$. Indeed, denoting $|\b|=l$ and estimating the
coefficients
\begin{equation*}
 \tex{
    \big|\sum\limits_{|\b|=l} \frac{(-1)^{l}}{\b!}  D^\b F(y) z_1^{\b_1}\cdots z_{N}^{\b_{N}}\big| \leq b_l
    |z|^l,
    }
\end{equation*}
by Stirling's formula we have that, for $l\gg 1$,
\begin{equation*}
%%\label{s11}
 \tex{
    b_l = \frac{N^l}{l!} \sup_{y\in \ren,
    |\b|=l} |D^\b F(y)| \approx \frac{N^l}{l!} l^{-l/10}
    {\mathrm e}^{l/10} \approx l^{-9l/10}c^l ={\mathrm e}^{-l \ln 9l/10 +l \ln c}.
    }
\end{equation*}
Note that, the series
%%\begin{equation*}
 %%\tex{
    $\sum b_l |z|^l$
   %% }
%%\end{equation*}
has the radius of convergence  $R=\infty$.
 Thus, we obtain the
following representation of the solution:
\begin{equation*}
%%\label{s12}
 \tex{
    w(y,\tau)= \sum\limits_{(\b)}  {\mathrm e}^{-\frac{|\b|}{10}\, \t} M_\b(u_0) \psi_\b(y),
    \quad \mbox{where} \quad
    \l_\b := -\frac{|\b|}{10}
    }
\end{equation*}
 and $\{\psi_\b\}$ are the eigenvalues and
eigenfunctions of the operator ${\bf B}$, respectively, and
\begin{equation*}
 \tex{
    M_\b(u_0):= \frac{1}{\sqrt{\b!}} \int\limits_{\ren} z_1^{\b_1}
    \cdots z_{N}^{\b_{N}} u_0(z) \, {\mathrm d}z
 }
\end{equation*}
are the corresponding momenta of the initial datum $w_0$ defined
by \eqref{s8}.

%%%%%%%%%%%%%%%%%%%%%%%%%%%%%%%%%%%%%%%%%%%%%
\subsection{Main spectral properties of the pair $\{\BB,\,
\BB^*\}$}

 Thus,  the
following holds \cite{EGKP}:

\begin{theorem}
\label{Th s1} {\rm (i)} The spectrum of ${\bf B}$ comprises real
eigenvalues only with the form
\begin{equation*}
%%\label{s13}
 \tex{
    \s({\bf B}):=\big\{\l_\b = -\frac{|\b|}{10}\,,\,|\b|=0,1,2,...\big\}.
    }
\end{equation*}
Eigenvalues $\l_\b$ have finite multiplicity with eigenfunctions,
\begin{equation}
\label{s14}
 \tex{
    \psi_\b(y):= \frac{(-1)^{|\b|}}{\sqrt{\b!}} D^\b F(y) \equiv \frac{(-1)^{|\b|}}{\sqrt{\b!}}
    \big(\frac{\p}{\p y_1}\big)^{\b_1}\cdots \big(\frac{\p}{\p y_N}\big)^{\b_N} F(y).
    }
\end{equation}

\noi{\rm (ii)} The subset of eigenfunctions
  $  \Phi=\{\psi_\b\}$
is complete in $L^2(\ren)$ and in $L_{\rho}^2(\ren)$.

\noi{\rm (iii)} For any $\l \notin \s({\bf B})$, the resolvent
  $  ({\bf B}-\l I)^{-1}$
is a compact operator in $L_{\rho}^2(\ren)$.
\end{theorem}

Subsequently, it was also shown in \cite{EGKP} that the adjoint
(in the dual metric of $L^2(\ren)$) operator of ${\bf B}$ given by
\begin{equation*}
 %%\label{B*}
 \tex{
    {\bf B}^*:=\D^5 - \frac{1}{10}\,\, y \cdot \nabla,
    }
\end{equation*}
 in the weighted space
$L_{{\rho}^*}^2(\ren)$, with the exponentially decaying weight
function
\begin{equation*}
%% \label{rho*}
 \tex{
    {\rho}^*(y) \equiv \frac{1}{\rho(y)} = {\mathrm e}^{-a|y|^{10/9}} >0,
    }
\end{equation*}
is a bounded linear operator,
\begin{equation*}
    {\bf B}^*:H_{{\rho}^*}^{10}(\ren) \to L_{{\rho}^*}^2(\ren),
    \,\, \mbox{so} \,\,
    \big\langle {\bf B} v, w\big\rangle = \big\langle v, {\bf B}^* w\big\rangle,
    \,\,
    v \in H_{\rho}^{10}(\ren), \,\,
    w \in H_{{\rho}^*}^{10}(\ren).
\end{equation*}
Moreover, the following theorem establishes the spectral properties of the
adjoint operator which will be very similar to those ones shown in Theorem\,\ref{Th s1}
for the operator $\bf{B}$.

\begin{theorem}
\label{Th s2}
{\rm (i)} The spectrum of ${\bf B}^*$ consists of eigenvalues of
finite multiplicity,
\begin{equation*}
%%\label{s15}
 \tex{
    \s({\bf B}^*)=\s({\bf B}):=\big\{\l_\b = -\frac{|\b|}{10}\,,\,|\b|=0,1,2,...\big\},
    }
\end{equation*}
and the eigenfunctions
   $\psi_\b^*(y)$
are polynomials of order $|\b|$.

\noi{\rm (ii)}  The subset of eigenfunctions
$    \Phi^*=\{\psi_\b^*\}$
is complete in $L_{{\rho}^*}^2(\ren)$.

\noi{\rm (iii)} For any $\l \notin \s({\bf B}^*)$, the resolvent
  $  ({\bf B}^*-\l I)^{-1}$
is a compact operator in $L_{{\rho}^*}^2(\ren)$.
\end{theorem}

It should be pointed out that, since $\psi_0=F$ and $\psi_0^*
\equiv 1$, we have
\begin{equation*}
 \tex{
     \langle \psi_0, \psi_0^* \rangle= \int\limits_{\ren} \psi_0\, {\mathrm
     d}y
      =\int\limits_{\ren} F (y) \, {\mathrm
     d}y=1.
     }
\end{equation*}
However, thanks to \eqref{s14}, we have that
\begin{equation*}
 \tex{
     \int\limits_{\ren} \psi_\b \equiv \langle \psi_\b, \psi_0^* \rangle =0 
     \quad \hbox{for any} \quad
     |\b|\neq 0.
     }
\end{equation*}
This expresses the orthogonality property to the adjoint
eigenfunctions in terms of the dual inner product.
 %% Due to
 %%Theorem\,\ref{Th s2}, the adjoint eigenfunctions are polynomials
 %%which form a complete subset in $L_{{\rho}^*}^2(\ren)$ with
 %%exponential decaying weight ${\rho}^*(y) = {\mathrm
 %%e}^{-a|y|^{4/3}}$.
 %%\par

Note that \cite{EGKP}, for the eigenfunctions $\{\psi_\b\}$ of
$\bf{B}$ denoted by \eqref{s14}, the corresponding adjoint
eigenfunctions are {\em generalized Hermite  polynomials} given by
\begin{equation}
\label{s16}
 \tex{
    \psi_\b^*(y):=\frac{1}{\sqrt{\b!}}\Big[y^\b + \sum\limits_{j=1}^{[|\b|/10]}
    \frac{1}{j!}\, \D^{5j} y^\b\Big].
    }
\end{equation}
Hence, the orthonormality condition holds
\begin{equation*}
%%\label{s17}
    \big\langle \psi_\b,\psi_\g \big\rangle=\d_{\b\g}
   \quad \hbox{for any} \quad \b,\;\g,
\end{equation*}
where $\big\langle \cdot,\cdot \big\rangle$ is the duality product
in $L^2(\ren)$ and $\d_{\b\g}$ is  Kronecker's delta. Also,
operators $\bf{B}$ and $\bf{B}^*$ have zero Morse index (no
eigenvalues with positive real parts are available).
 Key spectral results  can be extended \cite{EGKP} to
$2m$th-order linear poly-harmonic flows
 \begin{equation*}
%%  \label{PP1}
  \tex{
 u_t= - (-\D)^m u \quad \mbox{in} \quad \ren \times \re_+,
  }
  \end{equation*}
  where
 the elliptic equation for the rescaled kernel $F(y)$ takes the
 form
\begin{equation*}
%%\label{s5}
  \tex{
    {\bf B} F \equiv -(-\D_y)^m F + \frac{1}{2m}\,y \cdot \nabla_y F  +\frac{N}{2m} \,F=0
    \quad \hbox{in} \quad \ren,\quad \int\limits_{\ren} F(y) \, {\mathrm
    d}y=1.
     }
\end{equation*}
In particular, for $m=1$, we find the \emph{Hermite operator} and
the {\em Gaussian kernel} (see \cite{BS} for further information)
\begin{equation*}
 \tex{
    {\bf B} F \equiv \D F + \frac{1}{2}\, y\cdot \n F  +\frac{N}{2} \,F=0
     \LongA F(y)= \frac 1{(4 \pi)^{N/2}} \, {\mathrm
     e}^{-\frac{|y|^2}4},
    }
\end{equation*}
whose name is connected to fundamental works of Charles Hermite
on orthogonal polynomials $\{H_\b\}$ about 1870.
 These classic Hermite polynomials are obtained by differentiating
 the Gaussian: up to normalization constants,
 \be
 \label{Ga1}
  \tex{
  D^\b{\mathrm
     e}^{-\frac{|y|^2}4}= H_\b(y) \, {\mathrm
     e}^{-\frac{|y|^2}4} \quad \mbox{for any} \,\,\, \b.
     }
     \ee
 Note that, for
$N=1$, such operators and polynomial eigenfunctions in 1D were
studied earlier by Jacques C.F.~Sturm in 1836; on this history and
Sturm's main original calculations, see \cite[Ch.~1]{GalGeom}.

The generating formula \ef{Ga1} for (generalized) Hermite
polynomials is not available if $m \ge 2$, so that \ef{s16} are
obtained via a different procedure, \cite{EGKP}.

%%%%%%%%%%%%%%%%%%%%%%%%%%%%%%%%%%%%%%%%%%%%%%%%%%%%%
\setcounter{equation}{0}
\section{Similarity profiles for the Cauchy problem via  $n$-branching}
\label{S4}

%%%%%%%%%%%%%%%%%%%%%%%%%%%%%%%%%%%%%%%%%%%%%%%%

 In general, the construction of oscillatory similarity solutions of the
Cauchy problem for the TFE--10 \ef{i1} is a  difficult nonlinear
problem, which is harder than for the corresponding  FBP one.

On the other hand, for $n=0$, such similarity profiles  exist and are
given by eigenfunctions $\{\psi_\b\}$. In particular, the first
mass-preserving profile is just the rescaled kernel $F(y)$, so it
is unique, as was shown in Section \ref{S3}.

Hence, somehow, a possibility to visualize such an oscillatory
first ``nonlinear eigenfunction" $f(y)$ of changing sign, which
satisfies the {\em nonlinear eigenvalue problem} \eqref{self1},
  at least, for
sufficiently small $n > 0$ can be expected.

This suggests that, via an $n$-branching
approach argument, it is possible to ``connect" $f$ with the rescaled fundamental
profile $F$, satisfying the corresponding linear equation \ef{s4},
with all the necessary properties of $F$  presented in Section
\ref{S3}.

 Thus, we plan to describe the behaviour of the similarity
 profiles $\{f_\b\}$, as nonlinear eigenfunctions of \ef{self1}
for the TFE performing a ``homotopic" approach when $n\downarrow
0$ following a similar procedure performed in \cite{TFE4PabloVictor}.

 Homotopic approaches are well-known in the theory of
vector fields, degree, and nonlinear operator theory (see
\cite{Deim, KZ, VainbergTr} for details). However, we shall be
less precise in order to apply that approach, and here, a
``homotopic path" just declares existence of a continuous
connection (a curve) of solutions $f \in C_0$ that ends up at
$n=0^+$ at the linear eigenfunction $\psi_0(y)=F(y)$ or further
eigenfunctions $\psi_\b(y) \sim D^\b F(y)$, as \ef{s14} claims.
\par

Using classical branching theory in the case of finite regularity
of nonlinear operators involved, we formally show that the
necessary orthogonality condition holds deriving the corresponding
{\em Lyapunov--Schmidt branching equation}.
We will try to be as
rigorous as possible in supporting the delivery of the
nonlinear eigenvalues $\{\a_k\}$.

It is worth
mentioning that  TFE theory  for \emph{free boundary
problems} (FBPs) with nonnegative solutions is well understood
nowadays (at least in 1D). The FBP setting assumes posing three
standard boundary conditions at the interface, and such a theory
has been developed in many papers since 1990. The mathematical
formalities and general setting of the CP is still not fully
developed and a number of problems still remain open. In fact, the
concept of proper solutions of the CP is still partially obscure,
and moreover it seems that any classic or standard notions of
weak-mild-generalized-... solutions fail in the CP setting.

  Various ideas  associated with extensions
of smooth order-preserving semigroups are well known to be
effective for second-order nonlinear parabolic PDEs, when such a
construction is naturally supported by the maximum principle. The
 analysis of higher-order equations such as \eqref{i1} is much
harder than the corresponding  second-order equations or those in
divergent form
\begin{equation}
 \label{PME4}
    u_t =-(|u|^{n} u)_{xxxx}\quad \hbox{in} \quad \re \times \re_+\,,
\end{equation}
(see \cite{GalRDE4n} for a countable branching of similarity
solutions for \ef{PME4})
 because of the lack of the maximum
principle, comparison, order-preserving, monotone,  and potential
properties of the quasilinear operators involved.

It is clear that the CP for the \emph{poly-harmonic equation of tenth-order} \eqref{lin5}
is well-posed and has a unique solution given by the convolution
 \begin{equation*}
%% \label{b11}
 u(x,t)=b(x-\cdot,t)\, * \, u_0(\cdot),
  \end{equation*}
   where $b(x,t)$ is the fundamental solution of the operator $D_t - \Delta^5$. By the
   apparent
    connection
between \eqref{i1} and \eqref{lin5} (when $n=0$), intuitively at
least,
 this analysis provides us with a way to understand the CP for the TFE-10 by
 using  the fact that the proper
 solution of the CP for \eqref{i1}, with the same initial data $u_0$, is that one which converges to the corresponding
unique solution of the CP for \eqref{s1}, as $n\rightarrow 0$.
Thus, we shall use the patterns occurring for $n=0$, as branching
points of nonlinear eigenfunctions, so some extra detailed
properties of this linear flow will be necessary.

Further extensions of solutions
for non-small $n>0$ require a novel essentially non-local
technique of such nonlinear analysis, which remains an open
problem.

%%%%%%%%%%%%%%%%%%%%%%%%%%%%%%%%%%%%%%%%%%%%%%%%%%%%%%%%%%%%%%%%%%%
\subsection{Nonlinear eigenvalues $\{\a_k\}$ and transversality
conditions for the nonlinear eigenfunctions $f_k$}

In this first part of the section we establish the conditions and terms necessary for the
expansions of the parameter $\a$ and the nonlinear eigenfunctions, as well as the transversality
oscillatory conditions for such nonlinear eigenfunctions.

This will allow us to obtain the desired countable number of
solutions \eqref{i1} for the similarity equation \eqref{self1} via
Lyapunov-Schmidt reduction through the subsequent analysis.

The nonlinear eigenvalues $\{\a_k\}$ are obtained according to
non-self-adjoint spectral theory from Section \ref{S3}. We then
use  the explicit expressions for the eigenvalues and
eigenfunctions of the linear eigenvalue problem \ef{LP1} given in
Theorem \ref{Th s1}, where we also need the main conclusions of
the ``adjoint" Theorem \ref{Th s2}.

  Thus, taking the corresponding
linear equation from \eqref{self1} with $n=0$, we find, at least,
formally, that
\begin{equation*}
 \tex{
  n=0: \quad  \mathcal{L}(\a)f:=\D^5 f +\frac{1}{10}\, y \cdot\nabla f  +\a f=0.
     }
\end{equation*}

Moreover, from that equation, combined with the eigenvalues expressions
obtained in the previous section, we ascertain the following
critical values for the parameter $\a_k=\a_k(n)$,
\begin{equation}
\label{bf4}
 \tex{
  n=0: \quad   \a_k(0) := -\l_k + \frac{N}{10} \equiv \frac{k+N}{10} \quad \hbox{for any} \quad k=0,1,2,\ldots,
     }
\end{equation}
where $\l_k$ are the eigenvalues  defined in Theorem\,\ref{Th s1},
so that
\begin{equation*}
 \tex{
    \a_0(0)= \frac{N}{10}, \; \a_1(0)= \frac{N+1}{10},\; \a_2(0)= \frac{N+2}{10},\ldots ,
    \a_k(0)= \frac{k+N}{10}\ldots \,.
     }
\end{equation*}
In particular, when $k=0$, we have that $\a_0(0)= \frac{N}{10}$ and
the eigenfunction satisfies
\begin{equation*}
    {\bf B} F=0, \quad \mbox{so that} \quad
    \ker \mathcal{L}(\a_0) = \mathrm{span\,}\{\psi_0\} \quad(\psi_0=F),
\end{equation*}
and, hence, since $\l_0=0$ is a simple eigenvalue for the operator
$\mathcal{L}(\a_0)= {\bf B}$, its algebraic multiplicity is 1. In
general, we find that
\begin{equation*}
%%\label{bf5}
 \tex{
    \ker\big({\bf B} + \frac{k}{10}\, I \big)= \mathrm{span\,}\{\psi_\b, \, |\b|=k
    \}, \quad \hbox{for any}
    \quad k=0,1,2,3,\cdots\,,
    }
\end{equation*}
where the operator ${\bf B} + \frac{k}{10}\, I$ is Fredholm of
index zero since it is a compact perturbation of the identity of
linear type with respect to $k$. In other words,
$R[\mathcal{L}(\a_k)]$ is a closed subspace of $L_{\rho}^2(\ren)$
and, for each $\a_k$,
\begin{equation*}
    \hbox{dim} \ker(\mathcal{L}(\a_k))< \infty \andA \hbox{codim}R[\mathcal{L}(\a_k)]< \infty.
\end{equation*}

Then, for small $n>0$ in \eqref{self1}, we can assume the
following asymptotic expansions
\begin{equation}
\label{br3}
    \a_k(n):= \a_k+ \mu_{1,k} n+ o(n),\quad \mbox{and}
 \ee
 \be
  \label{br3N}
     |f|^n \equiv  {\mathrm e}^{n\ln |f|}:=
     1 +n \ln |f|+o(n).
\end{equation}
 As customary in bifurcation-branching theory \cite{KZ, VainbergTr},
existence of an expansion such as \ef{br3} will allow one to get
further expansion coefficients in
\begin{equation*}
%%\label{br31}
    \a_k(n):= \a_k+ \mu_{1,k} n + \mu_{2,k} n^2+ \mu_{3,k} n^3 + ...\,,
 \end{equation*}
 as the regularity of nonlinearities allows and suggests, though the convergence
 of such an analytic series can be questionable and is not under
 scrutiny here.

Another principle question is that, for oscillatory sign changing
profiles $f(y)$, the last expansion \ef{br3N} cannot be understood
in the pointwise sense. However, it can be naturally expected to be valid
in other metrics such as weighted $L^2$ or Sobolev spaces, as in
Section \ref{S3}, that used to be appropriate for the functional
setting of the equivalent integral equation and for that with
$n=0$.

\vspace{0.4cm}

\noi\underline{\sc Transversality conditions}.
Let us explain  why a certain ``transversality" of zeros of
possible solutions $f(y)$ is of key importance. As we see the nonlinear operator in
\ef{self1} can be written in the following equivalent form:
\begin{equation}
\label{bf2}
     \tex{
   \D^5 f + \frac{1-\a n}{10}\,
    y \cdot \nabla f + \a f+\nabla \cdot ((|f|^n -1) \n \D^4 f)=0.}
\end{equation}
then, we have to use the expansion for small $n>0$
 \be
 \label{nn1}
 |f|^n-1 \equiv {\mathrm e}^{n \, \ln |f|}-1= 1+ n \,
 \ln|f|+...-1=n \, \ln|f|+...\,,
  \ee
  which is true pointwise on any set $\{|f| \ge \e_0\}$ for
  an arbitrarily small
 fixed constant $\e_0>0$. However, in a small neighbourhood of any
 {\em zero} of $f(y)$, the expansion \ef{nn1} is no longer true.
 Nevertheless, it remains true in a weak sense provided that this
 zero is sufficiently transversal in a natural sense, i.e.,
 \be
 \label{nn2}
  \tex{
  \frac{|f|^n-1}n \rightharpoonup \ln|f| \asA n \to 0^+
  }
  \ee
  in $L^\infty_{\rm loc}$, since then
 the singularity $\ln |f(y)|$ is not more than {\em logarithmic}
 and, hence, is locally integrable in
 \be
 \label{int1}
  \tex{
f=- \big( \D^5  + \frac{1-\a n}{10}\,
    y \cdot \nabla +( \a +a)I \big)^{-1} (\nabla \cdot ((|f|^n -1) \n \D^4 f)+a f),
    }
 \ee
  where $a>0$ is a parameter to be chosen so that the inverse
  operator (a resolvent value) is a compact one in a weighted space $L^2_\rho(\ren)$; see
  Section \ref{S3}. We will show therein that the
  spectrum of
  $$
   \tex{
   \mc{L}(\a,n):= \D^5  + \frac{1-\a n}{10}\,
    y \cdot \nabla +\a I,
    }
    $$
    is always discrete and, actually,
 \begin{equation*}
 %%\label{sp33}
  \tex{
 \s(\mc{L}(\a,n))=\big\{(1-\a n)\big(- \frac{k}{10}\big)+\a,
 \,k=0,1,2,...\big\},
 }
  \end{equation*}
 so that any choice of $a>0$ such that $a \not \in \s(\mc{L})$ is suitable in \ef{int1}.

 Equivalently we are dealing with the limit
 $$n \ln^2 |f| \rightharpoonup 0, \quad \hbox{as} \quad n \downarrow 0^+,$$
 at least in a very weak sense, since by the expansion \eqref{nn1} we have that
 $$
   \tex{
    \frac{|f|^n-1}n - \ln |f| = \frac 12 \, n \, \ln^2 |f|+...
    \, .
    }
    $$

   Note also that actually we deal, in \ef{int1}, with an easier
  expansion
   \be
   \label{nn3}
   (|f|^n-1) \n \D^4 f = (n \, \ln|f|+...) \n \D^4 f,
    \ee
    so that even if $f(y)$ does not vanish transversely at a
    zero surface, the extra multiplier $\n \D^4 f(y)$ in \ef{nn3},
    which is supposed to vanish as well, helps to improve the
    corresponding weak convergence.
Furthermore, it is seen from \ef{bf2} that, locally in space
variables, the operator in \ef{int1} (with $a=0$ for simplicity)
acts like a standard Hammerstein--Uryson compact integral operator
with a sufficiently smooth kernel:
 \be
 \label{nn4}
 f \sim (\n \D^4)^{-1}[(|f|^n-1) \n \D^4 f].
  \ee

 Therefore, in order to justify our asymptotic branching analysis,
 one needs in fact to introduce such a functional setting and a class of solutions
 $$
 \mathcal P=\{f=f(\cdot,n): \,\,
  f \in H^{10}_{\rho}(\ren)\}, \quad \mbox{for which}
  $$
  %%for which:
 \be
 \label{nn5}
  \tex{
  {\mathcal P}: \quad \,\,
   (\n \D^4)^{-1}\big( \frac{|f|^n-1}n \,\n \D^4 f\big)  \to (\n \D^4)^{-1}(\ln \, |f|
 \n \D^4 f)
    \asA n
   \to 0^+
   }
   \ee
a.e. This is the precise statement on the regularity of possible solutions,
which is necessary to perform our asymptotic branching analysis.
 In 1D or in the radial
 geometry in $\ren$, \ef{nn5} looks rather constructive. However, in general, for
 complicated solutions with unknown types of compact supports in $\ren$,
 functional settings that can guarantee \ef{nn5}
 are not achievable still.
We mention again that, in particular, our formal analysis aims to
establish structures of difficult multiple zeros of the nonlinear
eigenfunctions $f_\g(y)$, at which \ef{nn5} can be violated, but
hopefully not in the a.e. sense.

Then, since \ef{br3N} is obviously pointwise violated at the
nodal set $\{f=0\}$ of $f(y)$, this imposes some restrictions on
the behaviour   of corresponding eigenfunctions $\psi_\b(y)$
($n=0$) close to their zero sets.
 Using well-known asymptotic and other related properties of the
{\em radial} analytic rescaled kernel $F(y)$ of the fundamental
solutions \ef{s3}, the generating formula of eigenfunctions
\ef{s14} confirms that the nodal set of analytic eigenfunctions
$\{\psi_\b=0\}$
 consists of isolated zero surfaces, which are ``transversal", at least in the a.e. sense,
with the only accumulation point at $y= \infty$.
 Overall, under such conditions, this indicates that
  \be
  \label{log1}
 \mbox{expansion (\ref{br3N}) contains not more than
 ``logarithmic" singularities a.e.},
  \ee
   which well suited the integral compact operators involved in
   the branching analysis.
  %%though we are far from claiming this as any rigorous issue.

 Moreover, when $n>0$ is not small enough, such an analogy and statements like
  \ef{log1} become
unclear, and global extensions of continuous $n$-branches induced by
some compact integral operators, i.e., nonexistence of turning
(saddle-node) points in $n$, require, as usual, some unknown
monotonicity-like results.

Then, in order to carry out our homotopic approach
we assume the expansion \ef{br3N} away from possible zero
surfaces of $f(y)$, which, by transversality, can be localized in
arbitrarily small neighbourhoods.

Indeed, it is clear that when $$|f| > \d >0, \quad \hbox{for any}
\quad \d>0,$$ there is no problem in approximating $|f|^n$ by
\eqref{br3N}, i.e.,
 $$ |f|^n =1+ O(n) \quad \hbox{as} \quad  n
\rightarrow 0^+.
 $$
 However, when
  $$
  |f| \leq \d, \quad \hbox{for any} \quad \d>0,
    $$
sufficiently small, the proof of such an approximation in weak
topology (as suffices for dealing with equivalent integral
equations) is far from clear unless
\begin{equation*}
 \mbox{the zeros of the $f$'s are also transversal a.e.},
  \end{equation*}
with a standard accumulating property at the
only interface zero surface. The latter issues have been studied
and described in \cite{EGK2} in the radial setting.
 Hence, we can suppose that such nonlinear eigenfunctions $f(y)$ are oscillatory
and infinitely sign changing close to the interface surface.

Therefore, if we assume that their zero surface is transversal
a.e. with a known geometric-like accumulation at the interface, we
find that, for any $n$ close to zero and any $\d= \d(n)
>0$ sufficiently small,
\begin{equation*}
    n| \ln |f| | \gg 1, \quad \hbox{if} \quad |f| \leq \d(n),
\end{equation*}
and, hence, on such subsets, $f(y)$ must be exponentially small:
\begin{equation*}
 \tex{
    | \ln |f| | \gg \frac{1}{n}\; \Longrightarrow \;\ln |f| \ll -\frac{1}{n}\;
    \Longrightarrow \; |f|  \ll {\mathrm e}^{-\frac{1}{n}}.
    }
\end{equation*}

Recall that this happens in  also exponentially small
neighbourhoods of the transversal zero surfaces.

Overall, using the periodic structure of the oscillatory component
at the interface \cite{EGK2} (we must admit that such delicate
properties of oscillatory structures of solutions are known for
the 1D and radial cases only, though we expect that these
phenomena are generic), we can control the singular coefficients
in (\ref{br3N}), and, in particular, to see that
\begin{equation}
 \label{f1loc}
    \ln |f| \in L^1_{\rm loc} (\ren).
\end{equation}
However, for most general geometric configurations of nonlinear
eigenfunctions $f(y)$, we do not have a proper proof of
(\ref{f1loc}) or similar estimates, so our further analysis is
still essentially formal.

\subsection{Derivation of the branching equation}

Under the above-mentioned transversality conditions and assuming the expansions
\eqref{br3}, for the nonlinear eigenvalues $\a_k$, and \eqref{br3N}, for the nonlinear
eigenfunctions f, we are able to obtain the branching equation applying the
classical Lyapunov-Schmidt method.

It is worth recalling again that our
computations below are to be understood as those dealing with the
equivalent integral equations and operators, so, in particular, we
can use the powerful facts on compactness of the resolvent
$(\BB-\l I)^{-1}$ and of the adjoint one $(\BB^*-\l I)^{-1}$ in
the corresponding weighted $L^2$-spaces. Note that, in such an
equivalent integral representation, the singular term in
(\ref{br3N}) satisfying (\ref{f1loc}) makes no principal
difficulty, so the  expansion  (\ref{br3N}) makes rather usual
sense for applying standard nonlinear operator theory.

 Thus, under natural assumptions, substituting \eqref{br3} into \eqref{self1}, for any
$k=0,1,2,3,\cdots$\,, we find that, omitting $o(n)$ terms when
necessary,
\begin{equation*}
 \tex{
    \nabla \cdot[(1 +n \ln |f|) \n \D^4 f]
    +\frac{1-\a_k n - \mu_{1,k} n^2}{10}\, y \cdot\nabla  f
    +(\a_k+ \mu_{1,k} n) f=0\,,
    }
\end{equation*}
and, rearranging terms,
\begin{equation*}
 \tex{
    \D^5 f +n\nabla \cdot( \ln |f| \n \D^4 f)
    +\frac{1}{10}\, y \cdot\nabla  f
    -\frac{\a_k n + \mu_{1,k} n^2}{10}\, y \cdot\nabla  f
    +\a_k f+ \mu_{1,k} n f=0\,.
    }
\end{equation*}
Hence, we finally have
\begin{equation*}
%%\label{br4}
 \tex{
    \big({\bf B} + \frac{k}{10}\,I\big) f +n \big[\nabla \cdot( \ln |f| \n \D^4 f)
    -\frac{\a_k}{10}\, y \cdot\nabla  f
    + \mu_{1,k} f \big]+o(n)   =0\,,
    }
\end{equation*}
which can be written in the following form:
\begin{equation}
\label{br5}
 \tex{
    \big({\bf B} + \frac{k}{10}\,I \big) f +n \mathcal{N}_k (f)+o(n)   =0\,,
     }
\end{equation}
with the operator
\begin{equation*}
%%\label{br6}
 \tex{
    \mathcal{N}_k (f):=\nabla \cdot( \ln |f| \n \D^4 f)
    -\frac{\a_k}{10}\, y \cdot\nabla  f
    + \mu_{1,k} f\,.
    }
\end{equation*}
Subsequently, as was shown in Section \ref{S3}, we have that
\begin{equation*}
%%\label{br7}
 \tex{
    \ker\big({\bf B} + \frac{k}{10}\, I \big)= \mathrm{span\,}\{\psi_\b\}_{|\b|=k}\quad \hbox{for any}
    \quad k=0,1,2,3,\cdots,
     }
\end{equation*}
where the operator
   $ {\bf B} + \frac{k}{10}\, I$
is Fredholm of index zero and
\begin{equation*}
 \tex{
    \dim \ker\big({\bf B} + \frac{k}{10}\, I\big)=  M_k  \ge 1 \quad \hbox{for any}
    \quad k=0,1,2,3,\cdots,
     }
\end{equation*}
where $M_k$ stands for the length of the vector $\{D^\b v, \,
|\b|=k\}$, so that $M_k>1$ for  $k \ge 1$.

Subsequently, we shall compute the coefficients involved in the
expansions \eqref{br3} and \eqref{br3N} applying the classical Lyapunov--Schmidt
method to \eqref{br5}  (branching approach when $n\downarrow 0$),
and, hence, describing the behaviour of the global solutions for
at least small values of the parameter $n>0$. Two cases are
distinguished. The first one in which the eigenvalue is simple and
the second for which the eigenvalues are semisimple. Note that due
to Theorems\;\ref{Th s1} and \ref{Th s2}, for any $k\geq 0$, the
algebraic multiplicities are equal to the geometric ones, so we do
not deal with the problem of introducing the generalized
eigenfunctions (no Jordan blocks are necessary for restrictions to
eigenspaces).

\ssk

\noi\underline{\sc Simple eigenvalue for $k=0$}. Since 0 is a
simple eigenvalue of ${\bf B}$ when $k=0$, i.e.,
\begin{equation*}
    \ker\,{\bf B} \oplus R[{\bf B}] = L_{\rho}^2(\ren),
\end{equation*}
the study of the case $k=0$ seems to be simpler than for other
different $k$'s because the dimension of the eigenspace is
$M_0=1$.

Thus, we shall describe the behaviour of solutions for
small $n>0$ and apply the classical Lyapunov--Schmidt method to
\eqref{br5} (assuming, as usual, some extra necessary regularity), in order to accomplish the
branching approach as $n\downarrow 0$, in two steps, when $k=0$
and $k$ is different from $0$.

Thus, owing to Section \ref{S3}, we already know that $0$ is
a simple eigenvalue of ${\bf B}$, i.e., $\ker\,{\bf B}=
\mathrm{span\,}\{\psi_0\}$ is one-dimensional. Hence, denoting by
$Y_0$ the complementary invariant subspace, orthogonal to
$\psi_0^*$, we set
\begin{equation*}
    f=\psi_0+V_0,
\end{equation*}
where $V_0 \in Y_0$.

Moreover, according to the spectral
properties of the operator ${\bf B}$, we define $P_0$ and $P_1$
such that $P_0+P_1=I$, to be the projections onto $\ker\,{\bf B}$
and $Y_0$ respectively. Finally, setting
\begin{equation}
\label{br8}
 V_0:=n \Phi_{1,0} + o(n),
\end{equation}
substituting the expression for $f$ into \eqref{br5} and passing
to the limit as $n\rightarrow 0^+$ leads to a linear inhomogeneous
equation for $\Phi_{1,0}$,
\begin{equation}
\label{br9}
    {\bf B}\Phi_{1,0}=- \mathcal{N}_0 (\psi_0),
\end{equation}
since ${\bf B} \psi_0=0$.

Furthermore, by  Fredholm theory, $V_0 \in
Y_0$ exists if and only if the right-hand side is orthogonal to
the one dimensional kernel of the adjoint operator ${\bf B}^*$
with $\psi_0^*=1$, because of \eqref{s16}. Hence, in the topology
of the dual space $L^2$, this requires the standard orthogonality
condition:
\begin{equation}
\label{br10}
    \big\langle \mathcal{N}_0 (\psi_0), 1\big\rangle=0.
\end{equation}
Then, \eqref{br9} has a unique solution $\Phi_{1,0} \in Y_0$
determining by \eqref{br8} a bifurcation branch for small $n>0$.
In fact, the algebraic equation \ef{br10} yields  the following
explicit expression for the coefficient $\mu_{1,0}$ of the
expansion \ef{br3} for the first eigenvalue $\a_0(n)$:

\begin{equation*}
%%\label{br11}
    \mu_{1,0}
    :=
 \tex{
    \frac{\langle -\nabla \cdot( \ln |\psi_0| \n \D^4 \psi_0)
    +\frac{N}{100}\,  y \cdot\nabla  \psi_0,\psi_0^*\rangle}
    {\langle \psi_0,\psi_0^*\rangle}
 }
     =
 \tex{
    \langle -\nabla \cdot( \ln |\psi_0| \n \D^4 \psi_0)
    +\frac{N}{100}\,  y \cdot\nabla  \psi_0,\psi_0^*\rangle.
 }
\end{equation*}
Consequently, in the particular case of having simple eigenvalues we just obtain one branch of solutions
emanating at $n=0$.

\ssk

\noi\underline{\sc Multiple eigenvalues for $k \ge 1$}. Next we
ascertain the number of branches in the case when the eigenvalues
of the operator ${\bf B}$ are semisimple.

For any
$k\geq 1$, we  know that
\begin{equation*}
 \tex{
    \dim \ker\big({\bf B} + \frac{k}{10} \,I\big)=  M_k >1.
    }
\end{equation*}
  Hence, in order to perform a similar analysis to the one done for simple eigenvalues
  we have to use the full eigenspace expansion
\begin{equation}
\label{br12}
 \tex{
    f=\sum\limits_{|\b|=k} c_\b \hat{\psi}_\b +V_k,
 }
\end{equation}
for every $k\geq 1$. Currently, for convenience, we denote
$$
 \tex{
 \{\hat{\psi}_\b\}_{|\b|=k}=\{\hat \psi_1,...,\hat \psi_{M_k}\},
  }
   $$
the natural basis of the $M_k$-dimensional eigenspace
$\ker\big({\bf B} + \frac{k}{10}\, I\big)$ and set
 $$
  \tex{
  \psi_k =
\sum_{|\b|=k} c_\b \hat{\psi}_\b.
 }
 $$
 Moreover,
$$
 \tex{
 V_k \in Y_k \quad \hbox{and} \quad V_k=\sum_{|\b|>k} c_\b {\psi}_\b,
  }
  $$
where $Y_k$
is the complementary invariant subspace of $\ker\big({\bf B} +
\frac{k}{10}\, I\big)$.

Furthermore, in the same way, as we did for
the case $k=0$, we define the $P_{0,k}$ and $P_{1,k}$, for every
$k\geq 1$, to be the projections of $\ker\big({\bf B} +
\frac{k}{10}\, I\big)$ and $Y_k$ respectively. We also expand $V_k$
 as
\begin{equation}
\label{br13}
    V_k:=n \Phi_{1,k} + o(n).
\end{equation}
Subsequently, substituting \eqref{br12} into \eqref{br5} and
passing to the limit as $n\downarrow 0^+$, we obtain the following
equation:
\begin{equation}
\label{br14}
 \tex{
    \big({\bf B}+ \frac{k}{10}\,I\big)\Phi_{1,k}=- \mc{N}_k \big(\sum_{|\b|=k} c_\b {\psi}_\b\big),
    }
\end{equation}
under the natural ``normalizing" constraint
\begin{equation}
\label{br15}
 \tex{
    \sum\limits_{|\b|=k} c_\b=1 \quad (c_\b \ge 0).
    }
\end{equation}
Therefore, applying the Fredholm alternative, $V_k \in Y_k$ exists
if and only if the term on the right-hand side of \eqref{br14} is
orthogonal to $\ker\,\big({\bf B}+ \frac{k}{10}\,I\big)$.
Then, multiplying the right-hand side of \eqref{br14} by $\psi_\b^*$,
for every $|\b|=k$,  in the topology of the dual space $L^2$, we
obtain an algebraic system of $M_k+1$ equations and the same
number of unknowns, $\{c_\b, \, |\b|=k\}$ and $\mu_{1,k}$:
\be
 \label{alg1}
  \tex{
\big\langle \mc{N}_k (\sum_{|\b|=k} c_\b {\psi}_\b), \psi^*_\b
\big\rangle=0 \quad \mbox{for all} \quad  |\b|=k,
 }
 \ee
 which is indeed the Lyapunov--Schmidt branching equation
 \cite{VainbergTr}.
In general, such algebraic systems are assumed to allow us to
obtain the branching parameters and hence establish the number of
different solutions induced on the given $M_k$-dimensional
eigenspace as the kernel of the operator involved.

However, a full solution of the {\em non-variational} algebraic
system \eqref{alg1} is a very difficult issue, though we claim
that the number of branches is expected to be related to the
dimension of the eigenspace $\ker\,\big({\bf B}^*+
\frac{k}{10}\,I\big)$.

In order to obtain the number of possible branches and with  the
objective of avoiding excessive notation, we analyze two typical
cases.

\vspace{0.4cm}

%%%%%%%%%%%%%%%%%%%%%%%%%%%%%%%%%%%%%%%%%%
\ssk

\noi\underline{\sc Computations for branching of dipole solutions in 2D}

Firstly, we ascertain some expressions for those
coefficients in the case when $|\b|=1$, $N=2$, and $M_1=2$, so that, in our notations,
$\{\psi_\b\}_{|\b|=1}=\{\hat \psi_1, \hat \psi_2\}$.

Consequently,
in this case, we obtain the following algebraic system: the
expansion coefficients of $\psi_1=c_1 \hat \psi_1+c_2 \hat \psi_2$
satisfy
\begin{equation}
\label{br16}
    \left\{\begin{array}{l}
    c_1  \langle \hat \psi_1^*,h_1 \rangle- \frac{c_1 \a_1}{10}\,
     \langle \hat \psi_1^*,y \cdot\nabla  \hat{\psi}_1 \rangle +c_1 \mu_{1,1}
    +c_2  \langle \hat \psi_1^*,h_2 \rangle- \frac{c_2 \a_1}{10}\,
     \langle \hat \psi_1^*,y \cdot\nabla  \hat{\psi}_2 \rangle = 0,\\
    c_1  \langle \hat \psi_2^*,h_1 \rangle- \frac{c_1 \a_1}{10}\,
     \langle \hat \psi_2^*,y \cdot\nabla  \hat{\psi}_1 \rangle
    + c_2  \langle \hat \psi_2^*,h_2 \rangle- \frac{c_2 \a_1}{10}\,
     \langle \hat \psi_2^*,y \cdot\nabla  \hat{\psi}_2 \rangle+ c_2 \mu_{1,1}=0,\\
    c_1+c_2=1,
    \end{array}\right.
\end{equation}
where
\begin{equation*}
   h_1:= \nabla \cdot  [ \ln (c_1 \hat{\psi}_1+c_2\hat{\psi}_2) \n \D^4
   \hat{\psi}_1 ],\,\,
    h_2:= \nabla \cdot
     [ \ln (c_1 \hat{\psi}_1+c_2\hat{\psi}_2) \n \D^4 \hat{\psi}_2 ],
\end{equation*}
and, $c_1$, $c_2$, and $\mu_{1,1}$ are the coefficients that we
want to calculate, $\a_1$ is regarded as the value of the
parameter $\a$ denoted by \eqref{bf4} and dependent on the
eigenvalue $\l_1$, for which $\hat \psi_{1,2}$ are the associated
eigenfunctions, and $\hat \psi_{1,2}^*$ the corresponding adjoint
eigenfunctions. Hence, substituting the expression $c_2=1-c_1$
from the third equation into the other two, we have the following
nonlinear algebraic system
\begin{equation}
\label{br17}
    \left\{\begin{array}{l}
    0=N_1(c_1,\mu_{1,1})- c_1\frac{\a_1}{10} \,\big[
     \langle \hat \psi_1^*,y \cdot\nabla  \hat{\psi}_1 \rangle
    -  \langle \hat \psi_1^*,y \cdot\nabla  \hat{\psi}_2 \rangle\big], \ssk\\
    0=N_2(c_1,\mu_{1,1}) - c_1\frac{\a_1}{10}\,\big[
     \langle \hat \psi_2^*,y \cdot\nabla  \hat{\psi}_1 \rangle
    - \langle \hat \psi_2^*,y \cdot\nabla  \hat{\psi}_2 \rangle\big]+ \mu_{1,1},
    \end{array}\right.
\end{equation}
where
\begin{equation*}
  \begin{split} &
   \tex{
   N_1(c_1,\mu_{1,1}):= c_1  \langle \hat \psi_1^*,h_1 \rangle
   + \langle \hat \psi_1^*,h_2 \rangle
   - \frac{\a_1}{10}\,
     \langle \hat \psi_1^*,y \cdot\nabla  \hat{\psi}_2 \rangle
   -c_1  \langle \hat \psi_1^*,h_2 \rangle +c_1 \mu_{1,1},
 }
   \\ &
    \tex{
   N_2(c_1,\mu_{1,1}):=c_1  \langle \hat \psi_2^*,h_1 \rangle
   +  \langle \hat \psi_2^*,h_2 \rangle
   - \frac{\a_1}{10} \,\langle \hat \psi_2^*,y \cdot\nabla  \hat{\psi}_2 \rangle-
   c_1 \langle \hat \psi_2^*,h_2 \rangle-c_1 \mu_{1,1}
    }
   \end{split}
\end{equation*}
represent the nonlinear parts of the algebraic system, with $h_0$
and $h_1$  depending on $c_1$.
\par

Subsequently, to guarantee existence of solutions of the system
\eqref{br16}, we apply the Brouwer fixed point theorem to
\eqref{br17} by supposing that the values $c_1$ and $\mu_{1,1}$
are the unknowns, in a disc sufficiently big
$D_R(\hat{c}_1,\hat{\mu}_{1,1})$ centered in a possible
nondegenerate zero $(\hat{c}_1,\hat{\mu}_{1,1})$. Thus, we write
the system \eqref{br17} in the matrix  form
\begin{equation*}
    \binom{0}{0}= \left(\begin{array}{cc}
    -\frac{\a_1}{10}\,\big[
   \langle \hat \psi_1^*,y \cdot\nabla  \hat{\psi}_1 \rangle
  -  \langle \hat \psi_1^*,y \cdot\nabla  \hat{\psi}_2 \rangle\big] & 0\\
  -\frac{\a_1}{10} \,\big[ \langle \hat \psi_2^*,y \cdot\nabla  \hat{\psi}_1 \rangle
  - \langle \hat \psi_2^*,y \cdot\nabla  \hat{\psi}_2 \rangle\big]
    & 1\end{array}\right)\binom{c_1}{\mu_{1,1}} + \binom{N_1(c_1,\mu_{1,1})}
    {N_2(c_1,\mu_{1,1})}.
\end{equation*}
Hence, we have that the zeros of the operator
\begin{equation*}
    \mathcal{F}(c_1,\mu_{1,1}):= \mf{M} \binom{c_1}{\mu_{1,1}} + \binom{N_1(c_1,\mu_{1,1})}
    {N_2(c_1,\mu_{1,1})}
\end{equation*}
are the possible solutions of \eqref{br17}, where $\mf{M}$ is the
matrix corresponding to the linear part of the system, while
 $$
  (N_1(c_1,\mu_{1,1}),N_2(c_1,\mu_{1,1}))^T,
   $$
 corresponds to the
nonlinear part. The application $\mathcal{H}:\mathcal{A} \times
[0,1] \to \re$, defined by
\begin{equation*}
    \mathcal{H}(c_1,\mu_{1,1},t):= \mf{M}
    \binom{c_1}{\mu_{1,1}} + t\binom{N_1(c_1,\mu_{1,1})}
    {N_2(c_1,\mu_{1,1})},
\end{equation*}
provides us with a homotopy transformation from the function
$\mathcal{F}(c_1,\mu_{1,1})= \mathcal{H}(c_1,\mu_{1,1},1)$ to its
linearization
\begin{equation}
\label{br18}
    \mathcal{H}(c_1,\mu_{1,1},0):= \mf{M} \binom{c_1}{\mu_{1,1}}.
\end{equation}

Thus, the system \eqref{br17} possesses a nontrivial solution if
\eqref{br18} has a nondegenerate zero, in other words, if the next condition
is satisfied
\begin{equation}
\label{br19}
     \langle \hat \psi_1^*,y \cdot\nabla  \hat{\psi}_1 \rangle-
     \langle \hat \psi_1^*,y \cdot\nabla  \hat{\psi}_2 \rangle\neq 0.
\end{equation}
Note that, if the substitution would have been $c_1=1-c_2$, the
condition might also be
\begin{equation*}
     \langle \hat \psi_2^*,y \cdot\nabla  \hat{\psi}_2 \rangle-
     \langle \hat \psi_2^*,y \cdot\nabla  \hat{\psi}_1 \rangle\neq 0.
\end{equation*}
Then, under condition \eqref{br19}, the system \eqref{br17} can be
written in the form
\begin{equation*}
%%\label{br20}
    \binom{c_1-\hat{c}_1}{\mu_{1,1}-\hat{\mu}_{1,1}}=-\mathcal{M}^{-1}
    \binom{N_1(c_1,\mu_{1,1})-\hat{c}_1}
    {N_2(c_1,\mu_{1,1})- \hat{\mu}_{1,1}},
\end{equation*}
which can be interpreted as a fixed point equation. Moreover,
applying Brower's fixed point theorem, we have that
\begin{align*}
    \hbox{Ind}((\hat{c}_1,\hat{\mu}_{1,1}),\mathcal{H}(\cdot,\cdot,0)) & =
    \mathcal{Q}_{C_R (\hat{c}_1,\hat{\mu}_{1,1})}(\mathcal{H}(\cdot,\cdot,0))
    %%\\ &
    =\hbox{Deg}(\mathcal{H}(\cdot,\cdot,0),D_R (\hat{c}_1,\hat{\mu}_{1,1}))\\ &
    = \hbox{Deg}(\mathcal{F}(c_1,\mu_{1,1}), D_R (\hat{c}_1,\hat{\mu}_{1,1})),
\end{align*}
where $\mathcal{Q}_{C_R
(\hat{c}_1,\hat{\mu}_{1,1})}(\mathcal{H}(\cdot,\cdot,0))$ defines
the number of rotations of the function
$\mathcal{H}(\cdot,\cdot,0)$ around the curve $C_R
(\hat{c}_1,\hat{\mu}_{1,1})$ and
$\hbox{Deg}(\mathcal{H}(\cdot,\cdot,0),D_R
(\hat{c}_1,\hat{\mu}_{1,1}))$ denotes the topological degree of
$\mathcal{H}(\cdot,\cdot,0)$ in $D_R (\hat{c}_1,\hat{\mu}_{1,1})$.
Owing to classical topological methods, both are equal.
\par

Thus, once we have proved the existence of solutions, we achieve
some expressions for the coefficients required:
\begin{equation*}
    \left\{
    \begin{matrix}
    %%{l} \begin{array}{l}
    \mu_{1,1} =c_2 ( \langle \hat \psi_1^*+\hat \psi_2^*,h_1-h_2 \rangle
    - \frac{\a_1}{10}\,
     \langle \hat \psi_1^*+\psi_2^*
    ,y \cdot\nabla  \hat{\psi}_1-y \cdot\nabla  \hat{\psi}_2 \rangle)
    \ssk\\
    %%\hspace{0.8cm}
     -  \langle \hat \psi_1^*+\hat \psi_2^*,h_1 \rangle
    + \frac{\a_1}{10}\,
     \langle \hat \psi_1^*+\hat \psi_2^*
    ,y \cdot\nabla  \hat{\psi}_1  \rangle, \qquad\qquad\qquad\,\,
    %%% \end{array}
    \ssk \\
    %%\hspace{0.5cm}
    c_1
    =1-c_2.\,\,\,\qquad\qquad\qquad\qquad\qquad\qquad\qquad\qquad\qquad\qquad\qquad
      \end{matrix}
      \right.
\end{equation*}
The expressions for the coefficients in a general case might be
accomplished after some tedious calculations, otherwise similar to
those performed above.

Note that, in general, those nonlinear
finite-dimensional algebraic problems are rather complicated, and
the problem of an optimal estimate of the number of different
solutions remains open.

Moreover, reliable multiplicity results
are very difficult to obtain. We expect that this number should be
somehow related (and even sometimes coincides) with the dimension
of the corresponding eigenspace of the linear operators ${\bf
B}+\frac{k}{10}\, I$, for any $k=0,1,2,\ldots\,$. This is a
conjecture only, and may be too illusive; see further supportive
analysis presented below.
\par

 However, we devote the remainder of this section to a possible answer
to that conjecture, which is not totally complete though, since we
are imposing some conditions.
\par

Thus, in order to detect the number of solutions of the nonlinear
algebraic system \eqref{br16}, we proceed to reduce this system to
a single equation for one of the unknowns. As a first step,
integrating by parts in the terms in which $h_1$ and $h_2$ are
involved and rearranging terms in the first two equations of the
system \eqref{br16}, we arrive at
 {\small $$
 \left\{
\begin{matrix}
 \tex{
   -\int\limits_{\ren} \nabla \psi_1^* \cdot \ln (c_1 \hat{\psi}_1+c_2\hat{\psi}_2) \n \D^4
   (c_1 \hat{\psi}_1
 }
    \tex{
    + c_2 \hat{\psi}_2)
    } 
    %%\qquad\qquad\quad
   %% \\
 %%\tex{
 - c_1 \frac{\a_1}{10}\,
    \int\limits_{\ren} \hat \psi_1^* \, y \cdot\nabla  \hat{\psi}_1
    +c_1 \mu_{1,1}- c_2  \frac{\a_1}{10}\,
    \int\limits_{\ren} \hat\psi_1^* \, y \cdot\nabla  \hat{\psi}_2 = 0,
    %}     
    %%\ssk\ssk
    \\
 \tex{
   -\int\limits_{\ren} \nabla \hat \psi_2^*\cdot \ln (c_1 \hat{\psi}_1+c_2\hat{\psi}_2) \n \D^4
   (c_2 \hat{\psi}_1+ c_2 \hat{\psi}_2)
 } 
 %%\qquad\qquad\quad
   %%\\
 %%\tex{
     -c_1  \frac{\a_1}{10}\,
   \int\limits_{\ren} \hat \psi_2^*\, y \cdot\nabla  \hat{\psi}_1
    + c_2 \mu_{1,1}
 %}
      \tex{
    -c_2  \frac{\a_1}{10}\,
    \int\limits_{\ren} \hat\psi_2^* \, y \cdot\nabla  \hat{\psi}_2 =0.
     }
\end{matrix}
 \right.
 $$
 }
 By the third equation, we have that $c_1 = 1- c_2$, and   hence,
setting $$c_1 \hat{\psi}_1+c_2\hat{\psi}_2 =
\hat{\psi}_1+(\hat{\psi}_2- \hat{\psi}_1)c_2$$ and substituting
these into those new expressions for the first two equations of
the system, we find that
 {\small \be
 \label{br59}
 \left\{
\begin{matrix}
 \tex{
   -\int\limits_{\ren} \nabla \hat \psi_1^* \cdot
   \ln (\hat{\psi}_1+
 }
     \tex{
     (\hat{\psi}_2-\hat{\psi}_1)c_2) \n \D^4
   (\hat{\psi}_1+(\hat{\psi}_2- \hat{\psi}_1)c_2) +\mu_{1,1}
    -c_2 \mu_{1,1}
     }
    \\
 \tex{
     - \frac{\a_1}{10}\,
    \int\limits_{\ren} \hat \psi_1^* \, y \cdot\nabla  \hat{\psi}_1
    + c_2  \frac{\a_1}{10}\, \int\limits_{\ren} \hat{\psi}_1^*\,  y
    \cdot(\nabla  \hat{\psi}_1-\nabla \hat{\psi}_2)  = 0,
     } \qquad\qquad\quad\,\, \ssk\ssk
      \\
 \tex{
  - \int\limits_{\ren} \nabla \hat \psi_2^*\cdot
   \ln (\hat{\psi}_1+
   }
     \tex{
     (\hat{\psi}_2-\hat{\psi}_1)c_2) \n \D^4
   (\hat{\psi}_1+(\hat{\psi}_2- \hat{\psi}_1)c_2)
    + c_2 \mu_{1,1}
 }\qquad\,\,\,
     \\
 \tex{
     - \frac{\a_1}{10}\,
    \int\limits_{\ren} \hat \psi_2^* \, y \cdot\nabla  \hat{\psi}_1
    +c_2  \frac{\a_1}{10}\,
    \int\limits_{\ren} \hat \psi_2^*\, y
    \cdot(\nabla  \hat{\psi}_1-\nabla \hat{\psi}_2) =0.
    }\qquad\qquad\quad\,\,
\end{matrix}
 \right.
\end{equation}
}
Subsequently, adding both equations, we have that
\begin{align*}
   \mu_{1,1} &  =
 \tex{
     \int\limits_{\ren}  (\nabla \hat \psi_1^*+ \nabla \hat \psi_2^*) \cdot
   \ln (\hat{\psi}_1+  (\hat{\psi}_2-\hat{\psi}_1)c_2) \n \D^4
   (\hat{\psi}_1+(\hat{\psi}_2- \hat{\psi}_1)c_2)
   }
    \\ &
 \tex{
     + \frac{\a_1}{10}\,
    \int\limits_{\ren} (\psi_1^*+ \psi_2^*) \, y
    \cdot \nabla  \hat{\psi}_1
    - c_2  \frac{\a_1}{10}\,
    \int\limits_{\ren} (\hat \psi_1^*+ \hat \psi_2^*)\, y
    \cdot(\nabla  \hat{\psi}_2-\nabla \hat{\psi}_1).
    }
\end{align*}
Thus, substituting it into the second equation of \eqref{br59}, we
obtain the following equation with the single unknown $c_2$:
{\small \begin{align*}
%%\label{br60}
   &
    \tex{
    -c_2^2 \, \frac{\a_1}{10}\,
    \int\limits_{\ren} (\hat \psi_1^*+ \hat \psi_2^*)\,  y
    \cdot (\nabla  \hat{\psi}_2-\nabla \hat{\psi}_1) +
    c_2 \frac{\a_1}{10}(\,
    \int\limits_{\ren} (\hat\psi_1^*+ 2 \hat\psi_2^*) \, y
    \cdot \nabla  \hat{\psi}_1- \,
    \int\limits_{\ren} \hat\psi_2^*\,  y \cdot \nabla \hat{\psi}_2)
    }
    \\ &
 \tex{
     - \frac{\a_1}{10}\,
    \int\limits_{\ren} \hat\psi_2^*\,  y \cdot\nabla  \hat{\psi}_1
    +\int\limits_{\ren} \nabla \psi_2^*\cdot
   \ln (\hat{\psi}_1+ (\hat{\psi}_2-\hat{\psi}_1)c_2) \n \D^4
   (\hat{\psi}_1+(\hat{\psi}_2- \hat{\psi}_1)c_2)
   }
   \\ &
 \tex{
    +c_2 \int\limits_{\ren} (\nabla \hat \psi_1^* +\nabla \hat\psi_2^*)\cdot
   \ln (\hat{\psi}_1+ (\hat{\psi}_2-\hat{\psi}_1)c_2) \n \D^4
   (\hat{\psi}_1+(\hat{\psi}_2- \hat{\psi}_1)c_2) =0,
   }
\end{align*}
}
which can be written in the following way:
\begin{equation*}
%% \label{FF1}
    c_2^2 A + c_2 B + C + \o(c_2) \equiv \mf{F} (c_2) + \o
    (c_2)=0.
\end{equation*}
 Here, $\o(c_2)$ can be considered as a perturbation of the quadratic
form $\mf{F} (c_2)$ with the coefficients
defined by
{\small \begin{align*}
    &
     \tex{
     A:=  -\frac{\a_1}{10}\,
    \int\limits_{\ren} (\hat\psi_1^*+ \hat\psi_2^*) \, y
    \cdot (\nabla  \hat{\psi}_2-\nabla \hat{\psi}_1),
 }
    \\ &
    \tex{
  B:= \frac{\a_1}{10}(\,
    \int\limits_{\ren} (\hat\psi_1^*+ 2 \hat\psi_2^*) \, y
    \cdot \nabla  \hat{\psi}_1- \,
    \int\limits_{\ren} \hat\psi_2^* y \cdot \nabla \hat{\psi}_2),
    }
     \quad
  C:=
 \tex{
   - \frac{\a_1}{10}\,
    \int\limits_{\ren} \hat\psi_2^*\, y \cdot\nabla  \hat{\psi}_1,
    }
     \\ &
  \o(c_2):=
  \tex{
   \int\limits_{\ren} \nabla \hat\psi_2^*\cdot
   \ln (\hat{\psi}_1+ (\hat{\psi}_2-\hat{\psi}_1)c_2) \n \D^4
   (\hat{\psi}_1+(\hat{\psi}_2- \hat{\psi}_1)c_2)
   }
   \\ &
    \tex{
    +c_2 \int\limits_{\ren} (\nabla \hat\psi_1^*+\nabla \hat\psi_2^*)\cdot
   \ln (\hat{\psi}_1+ (\hat{\psi}_2-\hat{\psi}_1)c_2) \n \D^4
   (\hat{\psi}_1+(\hat{\psi}_2- \hat{\psi}_1)c_2).
   }
\end{align*}
}

Since, due to the normalizing constraint \eqref{br15}, $c_2 \in
[0,1]$, solving the quadratic equation $ \mf{F} (c_2)$ yields:
\begin{enumerate}
\item[(i)] $c_2=0 \Longrightarrow \mf{F}(0) = C $; \quad $(ii)$\; $c_2 =1 \Longrightarrow \mf{F}(1)= A+B+C $; and
\item[(iii)] differentiating $\mf{F}$ with respect to $c_2$, we obtain that
$\mf{F}'(c_2) = 2 c_2 A+B$. Then, the critical point of the
function $\mf{F}$ is $c_2^* = -\frac{B}{2A}$ and its image is $
\mf{F}(c_2^*)=- \frac{B}{4A}+C$.
\end{enumerate}

Consequently, the conditions that must be imposed in order to have
more than one solution (we already know the existence of at least
one solution) are as follows:
%%\begin{enumerate}
%%\item[
$$(a)\; C  (A+B+C) >0;\quad  
%%\item[
(b)\; C \big(-\frac{B}{4A}+C\big)<0;\quad \hbox{and}\quad
(c)\; 0<-\frac{B}{2A}<1.$$
%%\end{enumerate}
Note that, for $-\frac{B}{4A}+C= 0$, we have just a single
solution.
%%\par
Hence, considering the equation again in the form
 $
 \mf{F} (c_2) + \o
(c_2)=0
 $,
  where $\o (c_2)$ is a perturbation of the quadratic form
$\mf{F} (c_2)$, and bearing in mind that the objective is to
detect  the number of solutions of the system \eqref{br16}, we
need to control  somehow this perturbation.

Under the conditions
(a), (b), and (c), $\mf{F} (c_2)$ possesses exactly two solutions.
Therefore, controlling the possible oscillations of the
perturbation $\o (c_2)$ in such a way that
\begin{equation*}
    \left\| \o (c_2) \right\|_{L^\infty} \leq \mf{F}(c_2^*),
\end{equation*}
we can assure that the number of solutions for \eqref{br16} is
exactly two. This is  the dimension of the kernel of the operator
${\bf B}+\frac{1}{10}\, I$ (as we expected in our more general
conjecture).

 The above particular example shows how
difficult the questions on existence and multiplicity of
solutions for such non-variational branching problems are.

 Recall that the actual values of the coefficients $A$, $B$, $C$,
 and others, for which the number of solutions crucially depends on, are very
 difficult to estimate, even numerically, in view of the complicated
 nature of the eigenfunctions \ef{s14} involved. To say nothing of
 the nonlinear perturbation $\o(c_2)$.

\vspace{0.2cm}

%%%%%%%%%%%%%%%%%%%%%%%%%%%%%%%%%%%%%%%%%%%%%%%%%%%%%
%%%%%%%%%%%%%%%%%%%%%%%%%%%%%%%%%%%%%%%%%%
\ssk

\noi\underline{\sc Branching computations for $|\b|=2$}

\ssk

Overall, the above analysis provides us with some expressions for
the solutions for the self-similar equation \eqref{self1}
depending on the value of $k$. Actually, we can achieve those
expressions for every critical value $\a_k$, but again the
calculus gets rather difficult.

For the sake of completeness, we now analyze the
case $|\b|=2$ and $M_2=3$, so that $\{\psi_\b\}_{|\b|=2}=\{\hat
\psi_1, \hat \psi_2, \hat \psi_3\}$ stands for a basis of the
eigenspace $\ker\big({\bf B} + \frac{1}{5}\, I \big)$, with $k=2$
($\l_k=- \frac{k}{10}$).
Thus, in this case, performing in a similar way as was done for
\eqref{br16} with
$\psi_2= c_1 \hat \psi_1+c_2 \hat \psi_2+c_3 \hat \psi_3$,
we arrive at the following algebraic system:
{\small \begin{equation}
\label{br61}
    \left\{\begin{array}{l}
    \begin{array}{r}
     c_1  \langle \hat \psi_1^*,h_1 \rangle
    +c_2  \langle \hat \psi_1^*,h_2 \rangle
    +c_3  \langle \hat \psi_1^*,h_3 \rangle - \frac{c_1 \a_2}{10}\,
     \langle \hat \psi_1^*,y \cdot\nabla  \hat{\psi}_1 \rangle
    - \frac{c_2 \a_2}{10}\,
     \langle \hat \psi_1^*,y \cdot\nabla  \hat{\psi}_2 \rangle
    \\- \frac{c_3 \a_2}{10}\,
     \langle \hat \psi_1^*,y \cdot\nabla  \hat{\psi}_3 \rangle +c_1 \mu_{1,2}= 0,
    \ssk \\
     c_1  \langle \hat \psi_2^*,h_1 \rangle
    + c_2  \langle \hat \psi_2^*,h_2 \rangle
    + c_2  \langle \hat \psi_2^*,h_3 \rangle  - \frac{c_1
    \a_2}{10}\,
     \langle \hat \psi_2^*,y \cdot\nabla  \hat{\psi}_1 \rangle
    - \frac{c_2 \a_2}{10}\,
     \langle \hat \psi_2^*,y \cdot\nabla  \hat{\psi}_2 \rangle
     \\  - \frac{c_3 \a_2}{10}\,
     \langle \hat \psi_2^*,y \cdot\nabla  \hat{\psi}_3 \rangle+ c_2 \mu_{1,2}=0,
     \ssk \\
     c_1  \langle \hat \psi_3^*,h_1 \rangle
    + c_2  \langle \hat \psi_3^*,h_2 \rangle
    + c_2  \langle \hat \psi_3^*,h_3 \rangle - \frac{c_1
    \a_2}{10}\,
     \langle \hat \psi_3^*,y \cdot\nabla  \hat{\psi}_1 \rangle
    - \frac{c_2 \a_2}{10}\,
     \langle \hat \psi_3^*,y \cdot\nabla  \hat{\psi}_2 \rangle
     \\  - \frac{c_3 \a_2}{10}\,
     \langle \hat \psi_3^*,y \cdot\nabla  \hat{\psi}_3 \rangle+ c_3 \mu_{1,2}=0,
     \end{array}\\
    c_1+c_2+c_3=1,
    \end{array}\right.
\end{equation}}
where
 $$
  \begin{matrix}
   h_1:= \nabla \cdot  [ \ln (c_1 \hat{\psi}_1+c_2\hat{\psi}_2
   +c_3 \hat{\psi}_3) \n \D^4
   \hat{\psi}_1 ], %%\\
 \quad    h_2:= \nabla \cdot
     [ \ln (c_1 \hat{\psi}_1+c_2\hat{\psi}_2
     +c_3 \hat{\psi}_3) \n \D^4 \hat{\psi}_2 ], \\
   h_3:= \nabla \cdot  [ \ln (c_1 \hat{\psi}_1+c_2\hat{\psi}_2
   +c_3 \hat{\psi}_3) \n \D^4
   \hat{\psi}_3 ],
\end{matrix}
 $$
and $c_1$, $c_2$, $c_3$, and $\mu_{1,2}$ are the unknowns to be
evaluated. Moreover, $\a_2$ is regarded as the value of the parameter
$\a$ denoted by \eqref{bf4} and is dependent on the eigenvalue
$\l_2$ with $\hat{\psi}_1,\hat{\psi}_2,\hat{\psi}_3$ representing
the associated eigenfunctions and
$\hat{\psi}_1^*,\hat{\psi}_2^*,\hat{\psi}_3^*$ the corresponding
adjoint eigenfunctions.
\par
Subsequently, substituting $c_3=1-c_1-c_2$ into the first three
equations and performing an argument based upon the Brower fixed
point theorem and the topological degree as the one done above for
the case $|\b|=1$, we ascertain the existence of a nondegenerate
solution of the algebraic system if the following condition is
satisfied:
\begin{equation*}
%%\label{br62}
     \langle \hat \psi_1^*,y \cdot\nabla  (\hat{\psi}_3-\hat{\psi}_1) \rangle
     \langle \hat \psi_2^*,y \cdot\nabla  (\hat{\psi}_3-\hat{\psi}_2) \rangle-
     \langle \hat \psi_1^*,y \cdot\nabla  (\hat{\psi}_3-\hat{\psi}_2) \rangle
     \langle \hat \psi_2^*,y \cdot\nabla  (\hat{\psi}_3-\hat{\psi}_1) \rangle\neq 0.
\end{equation*}
Note that, by similar substitutions, other conditions might be
obtained.
\par
Furthermore, once we know the existence of at least one solution,
we proceed now with a possible way of computing the number of
solutions of the nonlinear algebraic system \eqref{br61}.
Obviously, since the dimension of the eigenspace is bigger than
that in the case $|\b|=1$, the difficulty in obtaining
multiplicity results increases.
\par
First, integrating by parts in the nonlinear terms, in which
$h_1$, $h_2$ and $h_3$ are involved, and rearranging terms in the
first three equations gives
{\small \begin{align*}
 \tex{
   -\int\limits_{\ren} \nabla \psi_1^* \cdot
   \ln (c_1 \hat{\psi}_1+c_2\hat{\psi}_2+ c_3 \hat{\psi}_3)}
   & \tex{ \n \D^4
   (c_1 \hat{\psi}_1+ c_2 \hat{\psi}_2+ c_3 \hat{\psi}_3)
 }
    \tex{
     - c_1 \frac{\a_2}{10}\,
    \int\limits_{\ren} \hat \psi_1^*\, y \cdot\nabla  \hat{\psi}_1+c_1 \mu_{1,2}
    }
    \\ &
 \tex{
    - c_2  \frac{\a_2}{10}\,
    \int\limits_{\ren} \hat\psi_1^*\, y \cdot\nabla  \hat{\psi}_2
    - c_3  \frac{\a_2}{10}\,
    \int\limits_{\ren} \hat\psi_1^*\, y \cdot\nabla  \hat{\psi}_3 =
    0,
    }
\end{align*}
\begin{align*}
 \tex{
   -\int\limits_{\ren} \nabla \hat \psi_2^*\cdot
   \ln (c_1 \hat{\psi}_1+c_2\hat{\psi}_2+ c_3 \hat{\psi}_3)} &
   \tex{ \n \D^4
   (c_2 \hat{\psi}_1+ c_2 \hat{\psi}_2+ c_3 \hat{\psi}_3)
 }
 \tex{
     -c_1  \frac{\a_2}{10}\,
   \int\limits_{\ren} \hat \psi_2^* \, y \cdot\nabla  \hat{\psi}_1
    + c_2 \mu_{1,2}
 }
     \\ &
      \tex{
    -c_2  \frac{\a_2}{10}\,
    \int\limits_{\ren} \hat\psi_2^*\,  y \cdot\nabla  \hat{\psi}_2
    -c_3  \frac{\a_2}{10}\,
   \int\limits_{\ren} \hat \psi_2^*\, y \cdot\nabla  \hat{\psi}_3=0,
     }
\end{align*}
\begin{align*}
 \tex{
   -\int\limits_{\ren} \nabla \hat \psi_3^*\cdot
   \ln (c_1 \hat{\psi}_1+c_2\hat{\psi}_2+ c_3 \hat{\psi}_3)} &
   \tex{ \n \D^4
   (c_2 \hat{\psi}_1+ c_2 \hat{\psi}_2+ c_3 \hat{\psi}_3)
 }
 \tex{
     -c_1  \frac{\a_2}{10}\,
   \int\limits_{\ren} \hat \psi_3^*\, y \cdot\nabla  \hat{\psi}_1
    + c_3 \mu_{1,2}
 }
     \\ &
      \tex{
    -c_2  \frac{\a_2}{10}\,
    \int\limits_{\ren} \hat\psi_3^*\, y \cdot\nabla  \hat{\psi}_2
    -c_3  \frac{\a_2}{10}\,
   \int\limits_{\ren} \hat \psi_3^* \, y \cdot\nabla  \hat{\psi}_3=0.
     }
\end{align*}}
According to the fourth equation, we have that $c_1 = 1-c_2-c_3$. Then,
setting
\begin{equation*}
    c_1 \hat{\psi}_1+c_2\hat{\psi}_2+ c_3 \hat{\psi}_3 =
    \hat{\psi}_1+c_2(\hat{\psi}_2-\hat{\psi}_1)+ c_3
    (\hat{\psi}_3-\hat{\psi}_1)
\end{equation*}
and substituting it into the expressions obtained above for the
first three equations of the system yield
{\small \begin{align*}
 \tex{
   -\int\limits_{\ren} }
    &
     \tex{
      \nabla \hat \psi_1^* \cdot
   \ln (\hat{\psi}_1+
     (\hat{\psi}_2-\hat{\psi}_1)c_2+(\hat{\psi}_3-\hat{\psi}_1)c_3) \n \D^4
   (\hat{\psi}_1+(\hat{\psi}_2- \hat{\psi}_1)c_2+
   (\hat{\psi}_3- \hat{\psi}_1)c_3)
     }
     \\ &
 \tex{
     +\mu_{1,2}
    -c_2 \mu_{1,2}- c_3 \mu_{1,2} - \frac{\a_2}{10}\,
    \int\limits_{\ren} \hat \psi_1^*\,  y \cdot\nabla  \hat{\psi}_1
    %% }
    %%\\ &
 %%\tex{
    + \frac{\a_2}{10} \,\int\limits_{\ren} \hat{\psi}_1^*\, y
    \cdot( (\nabla  \hat{\psi}_1-\nabla \hat{\psi}_2)c_2 +
    (\nabla \hat{\psi}_1 - \nabla\hat{\psi}_3)c_3)  = 0,
     }
\end{align*}
\begin{equation}
\begin{split}
\label{br63}
 \tex{
   -\int\limits_{\ren}  }
    &
     \tex{
      \nabla \hat \psi_2^* \cdot
   \ln (\hat{\psi}_1+
     (\hat{\psi}_2-\hat{\psi}_1)c_2+(\hat{\psi}_3-\hat{\psi}_1)c_3) \n \D^4
   (\hat{\psi}_1+(\hat{\psi}_2- \hat{\psi}_1)c_2+
   (\hat{\psi}_3- \hat{\psi}_1)c_3)
     }
     \\ &
 \tex{
     +c_2 \mu_{1,2}
    - \frac{\a_2}{10}\,
    \int\limits_{\ren} \hat \psi_2^*\,  y \cdot\nabla  \hat{\psi}_1
     %%}
    %%\\ &
 %%\tex{
    + \frac{\a_2}{10}\, \int\limits_{\ren} \hat{\psi}_2^*\, y
    \cdot( (\nabla  \hat{\psi}_1-\nabla \hat{\psi}_2)c_2 +
    (\nabla \hat{\psi}_1 - \nabla\hat{\psi}_3)c_3)  = 0,
     }
\end{split}
\ee
\begin{align*}
 \tex{
   -\int\limits_{\ren}  }
    &
     \tex{
      \nabla \hat \psi_3^* \cdot
   \ln (\hat{\psi}_1+
     (\hat{\psi}_2-\hat{\psi}_1)c_2+(\hat{\psi}_3-\hat{\psi}_1)c_3) \n \D^4
   (\hat{\psi}_1+(\hat{\psi}_2- \hat{\psi}_1)c_2+
   (\hat{\psi}_3- \hat{\psi}_1)c_3)
     }
     \\ &
 \tex{
     +c_3\mu_{1,2}
    - \frac{\a_2}{10}\,
    \int\limits_{\ren} \hat \psi_3^*\, y \cdot\nabla  \hat{\psi}_1
     %%}
    %%\\ &
 %%\tex{
    + \frac{\a_2}{10} \,\int\limits_{\ren} \hat{\psi}_3^*\,  y
    \cdot( (\nabla  \hat{\psi}_1-\nabla \hat{\psi}_2)c_2 +
    (\nabla \hat{\psi}_1 - \nabla\hat{\psi}_3)c_3)  = 0.
     }
\end{align*}}
Now, adding the first equation of \eqref{br63} to the other two,
we have that
{\small \begin{align*}
 \tex{
  - \int\limits_{\ren}   }
    &
     \tex{
     (\nabla \hat \psi_1^*+ \nabla \hat{\psi}_2^*) \cdot
   \ln (\hat{\psi}_1+
     (\hat{\psi}_2-\hat{\psi}_1)c_2+(\hat{\psi}_3-\hat{\psi}_1)c_3) \n \D^4
   (\hat{\psi}_1+(\hat{\psi}_2- \hat{\psi}_1)c_2+
   (\hat{\psi}_3- \hat{\psi}_1)c_3)
     }
     \\ &
 \tex{
     +\mu_{1,2}
    - c_3 \mu_{1,2} - \frac{\a_2}{10}\,
    \int\limits_{\ren} (\hat \psi_1^*+ \hat{\psi}_2^*) \, y \cdot\nabla  \hat{\psi}_1
    %% }
    %%\\ &
 %%\tex{
    + \frac{\a_2}{10}\, \int\limits_{\ren} (\hat \psi_1^*+ \hat{\psi}_2^*)\,  y
    \cdot( (\nabla  \hat{\psi}_1-\nabla \hat{\psi}_2)c_2 +
    (\nabla \hat{\psi}_1 - \nabla \hat{\psi}_3)c_3)  = 0,
     }
\end{align*}
\begin{align*}
 \tex{
  - \int\limits_{\ren} }
    &
     \tex{
      (\nabla \hat \psi_1^*+ \nabla \hat{\psi}_3^*) \cdot
   \ln (\hat{\psi}_1+
     (\hat{\psi}_2-\hat{\psi}_1)c_2+(\hat{\psi}_3-\hat{\psi}_1)c_3) \n \D^4
   (\hat{\psi}_1+(\hat{\psi}_2- \hat{\psi}_1)c_2+
   (\hat{\psi}_3- \hat{\psi}_1)c_3)
     }
     \\ &
 \tex{
     +\mu_{1,2}
    - c_2 \mu_{1,2} - \frac{\a_2}{10}\,
    \int\limits_{\ren} (\hat \psi_1^*+ \hat{\psi}_3^*)\, y \cdot\nabla  \hat{\psi}_1
    %% }
    %%\\ &
 %%\tex{
    + \frac{\a_2}{10}\, \int\limits_{\ren} (\hat \psi_1^*+ \hat{\psi}_3^*)\, y
    \cdot( (\nabla  \hat{\psi}_1-\nabla \hat{\psi}_2)c_2 +
    (\nabla \hat{\psi}_1 - \nabla \hat{\psi}_3)c_3)  = 0.
     }
\end{align*}
}
Subsequently, subtracting those equations yields
\begin{align*}
 \tex{
    \mu_{1,2}
 }
 &
 \tex{
     = \frac{1}{c_3-c_2} \,\big[ \int\limits_{\ren}
      (\nabla \hat \psi_3^*- \nabla \hat{\psi}_2^*) \cdot
   \ln \Psi \n \D^4 \Psi -\frac{\a_2}{10}\,
    \int\limits_{\ren} (\hat \psi_2^*- \hat{\psi}_3^*)\,  y \cdot\nabla  \hat{\psi}_1
  }
     \\ &
  \tex{
    +\frac{\a_2}{10}\, \int\limits_{\ren} (\hat \psi_2^*- \hat{\psi}_3^*)\,  y
    \cdot( (\nabla  \hat{\psi}_1-\nabla \hat{\psi}_2)c_2 +
    (\nabla \hat{\psi}_1 - \nabla \hat{\psi}_3)c_3) \big],
  }
\end{align*}
where $\Psi = \hat{\psi}_1+
(\hat{\psi}_2-\hat{\psi}_1)c_2+(\hat{\psi}_3-\hat{\psi}_1)c_3$.
Thus, substituting it into \eqref{br63} (note that, from the
substitution into one of the last two equations, we obtain the
same equation), we arrive at the following system, with $c_2$ and
$c_3$ as the unknowns:
{\small \begin{align*}
&  \tex{
   -c_3 \int\limits_{\ren} (\nabla \hat{\psi}_1^*-
   \nabla \hat{\psi}_2^* 
   %% }
    %%&
     %%\tex{
     +\nabla \hat{\psi}_3^*) \cdot
   \ln \Psi
     \n \D^4 \Psi + c_2 \int\limits_{\ren} (\nabla \hat{\psi}_1^*+
   \nabla \hat{\psi}_2^*-\nabla \hat{\psi}_3^*) \cdot
   \ln \Psi \n \D^4 \Psi} \\ & \tex{ + \int\limits_{\ren} (
   \nabla \hat{\psi}_3^*-\nabla \hat{\psi}_2^*) \cdot
   \ln \Psi  \n \D^4 \Psi - \frac{\a_2}{10} \int\limits_{\ren} (\hat{\psi}_2^*-\hat{\psi}_3^*)\,
   y \cdot\nabla \hat{\psi}_1 \,
   %%} \\ & \tex{
    +c_2 \frac{\a_2}{10}\,
    [\int\limits_{\ren} (\hat{\psi}_2^*-\hat{\psi}_3^*) \, y \cdot\nabla
    (2\hat{\psi}_1 -\hat{\psi}_2)+\int\limits_{\ren}
    \hat{\psi}_1^*\, y \cdot\nabla  \hat{\psi}_1]
     }
     \\ &
 \tex{
    + c_3 \frac{\a_2}{10}\, [\int\limits_{\ren} (\hat{\psi}_2^*-\hat{\psi}_3^*)\, y \cdot\nabla
    (2\hat{\psi}_1 -\hat{\psi}_3)- \int\limits_{\ren}
    \hat{\psi}_1^*\, y \cdot\nabla  \hat{\psi}_1]
   }
      %%\end{align*}
      %%\begin{align*}
 \\  & \tex{
     + c_2 c_3  \frac{\a_2}{10}[\, \int\limits_{\ren} \hat{\psi}_1^*\, y
    \cdot( \nabla  \hat{\psi}_3-\nabla \hat{\psi}_2) - \, \int\limits_{\ren}
    (\hat{\psi}_2^* - \hat{\psi}_3^*) \, y \cdot (2\nabla \hat{\psi}_1 -
    \nabla \hat{\psi}_2-\nabla \hat{\psi}_3)]} \\ &
    \tex{
    + c_3^2  \frac{\a_2}{10}\, \int\limits_{\ren} (\hat{\psi}_1^*-
    \hat{\psi}_2^*+\hat{\psi}_3^*) \, y
    \cdot( \nabla \hat{\psi}_1 - \nabla\hat{\psi}_3)
    %%}\\ &
    %%\tex{
    - c_2^2  \frac{\a_2}{10}\, \int\limits_{\ren} (\hat{\psi}_1^*+
    \hat{\psi}_2^*-\hat{\psi}_3^*)\, y
    \cdot(\nabla  \hat{\psi}_1-\nabla \hat{\psi}_2)     = 0,
    }
%%\end{align*}
%%\begin{align*}
 \\ & \tex{
   -c_3 \int\limits_{\ren} \nabla \hat{\psi}_2^* \cdot
   \ln \Psi
% }
    %&
     %%\tex{
     \n \D^4 \Psi +c_2 \int\limits_{\ren} \nabla \hat{\psi}_3^* \cdot
   \ln \Psi \n \D^4 \Psi - c_3 \frac{\a_2}{10}\,
    \int\limits_{\ren} \hat{\psi}_2^*\, y \cdot\nabla  \hat{\psi}_1
    +c_2 \frac{\a_2}{10}\,
    \int\limits_{\ren} \hat{\psi}_3^*\, y \cdot\nabla  \hat{\psi}_1
     }
     \\ &
 \tex{
    + c_3 \frac{\a_2}{10}\, \int\limits_{\ren} \hat{\psi}_2^*\, y
    \cdot( (\nabla  \hat{\psi}_1-\nabla \hat{\psi}_2)c_2 +
    (\nabla \hat{\psi}_1 - \nabla\hat{\psi}_3)c_3)
    %% }
     %%\\ &
 %%\tex{
    - c_2 \frac{\a_2}{10}\, \int\limits_{\ren} \hat{\psi}_3^*\, y
    \cdot( (\nabla  \hat{\psi}_1-\nabla \hat{\psi}_2)c_2 +
    (\nabla \hat{\psi}_1 - \nabla\hat{\psi}_3)c_3)  = 0.
     }
\end{align*}}
These can be re-written in the following form:
{\small \begin{equation}
\label{br64}
    \begin{split}
    &
    A_1 c_2^2+ B_1 c_3^2 + C_1 c_2 +D_1 c_3+ E_1 c_2 c_3 +\o_1 (c_2,c_3)=0,
    \\ &
    A_2 c_2^2+ B_2 c_3^2 + C_2 c_2 +D_2 c_3+ E_2 c_2 c_3 +\o_2 (c_2,c_3)=0,
    \end{split}
\end{equation}
 where
\begin{align*}
    \tex{
    \o_1 (c_2,c_3)} & \tex{
   :=- c_3 \int\limits_{\ren} (\nabla \hat{\psi}_1^*-
   \nabla \hat{\psi}_2^*
     +\nabla \hat{\psi}_3^*) \cdot
   \ln \Psi
     \n \D^4 \Psi
     %%}
     %%\\ & \tex{
      +c_2 \int\limits_{\ren} (\nabla \hat{\psi}_1^*+
   \nabla \hat{\psi}_2^*-\nabla \hat{\psi}_3^*) \cdot
   \ln \Psi \n \D^4 \Psi
   } \\ & \tex{
   + \int\limits_{\ren} (
   \nabla \hat{\psi}_2^*-\nabla \hat{\psi}_3^*) \cdot
   \ln \Psi \n \D^4 \Psi- \frac{\a_2}{10} \int\limits_{\ren} (\hat{\psi}_2^*-\hat{\psi}_3^*)\,
   y \cdot\nabla \hat{\psi}_1
   }
\end{align*}
%%\com{PAC: $-+$ two lines por arriva???}
  and
\begin{align*}
    \tex{
    \o_2 (c_2,c_3)} & \tex{
    :=-c_3 \int\limits_{\ren} \nabla \hat{\psi}_2^* \cdot
   \ln \Psi \nabla \D \Psi + c_2 \int\limits_{\ren} \nabla \hat{\psi}_3^* \cdot
   \ln \Psi \nabla \D \Psi
   }
\end{align*}}
are the perturbations of the quadratic polynomials
 $$
 \mf{F}_i(c_2,c_3) :=
A_i c_2^2+ B_i c_3^2 + C_i c_2 +D_i c_3+ E_i c_2 c_3,
\quad \hbox{with} \quad i=1,2.
 $$
  The coefficients
of those quadratic expressions are given by
{\small\begin{align*}
    &
     \tex{
     A_1:=  -\frac{\a_2}{10}\, \int\limits_{\ren} (\hat{\psi}_1^*+
    \hat{\psi}_2^*-\hat{\psi}_3^*) \,y
    \cdot(\nabla  \hat{\psi}_1-\nabla \hat{\psi}_2)  ,\quad 
 %%}
    %%\\ &
    %%\tex{
  B_1:= \frac{\a_2}{10}\, \int\limits_{\ren} (\hat{\psi}_1^*-
    \hat{\psi}_2^*+\hat{\psi}_3^*) \,y
    \cdot( \nabla \hat{\psi}_1 - \nabla\hat{\psi}_3) ,
    }
     \\ &
  C_1:=
 \tex{
   \frac{\a_2}{10}\,
    [\int\limits_{\ren} (\hat{\psi}_2^*-\hat{\psi}_3^*) \,y \cdot\nabla
    (2\hat{\psi}_1 -\hat{\psi}_2)+ \int\limits_{\ren}
    \hat{\psi}_1 y \cdot\nabla  \hat{\psi}_1] ,\quad 
    }
     %%\\ &
  D_1:=
  \tex{
   \frac{\a_2}{10}\, [\int\limits_{\ren} (\hat{\psi}_2^*-\hat{\psi}_3^*) \, y \cdot\nabla
    (2\hat{\psi}_1 -\hat{\psi}_3)- \int\limits_{\ren}
    \hat{\psi}_1\,  y \cdot\nabla  \hat{\psi}_1] ,
   }
   \\ &
    \tex{
   E_1:=
    \frac{\a_2}{10}[\, \int\limits_{\ren} \hat{\psi}_1^*\,  y
    \cdot( \nabla  \hat{\psi}_3-\nabla \hat{\psi}_2) - \, \int\limits_{\ren}
    ( \hat{\psi}_2^* - \hat{\psi}_3^*)\,  y \cdot (2\nabla \hat{\psi}_1 -
    \nabla \hat{\psi}_2-\nabla \hat{\psi}_3)],
   }
\end{align*}
\begin{align*}
    &
     \tex{
     A_2:=  -\frac{\a_2}{10}\, \int\limits_{\ren} \hat{\psi}_3^*\, y
    \cdot( \nabla  \hat{\psi}_1-\nabla \hat{\psi}_2),\quad
 %}
    %%\\ &
    %%\tex{
  B_2:= \frac{\a_2}{10}\, \int\limits_{\ren} \hat{\psi}_2^*\,  y
    \cdot(
    (\nabla \hat{\psi}_1 - \nabla\hat{\psi}_3) ,
    }
     \\ &
  C_2:=
 \tex{
   \frac{\a_2}{10}\,
    \int\limits_{\ren} \hat{\psi}_3^*\,  y \cdot\nabla  \hat{\psi}_1,\quad 
    }
     %%\\ &
  D_2:=
  \tex{
   -  \frac{\a_2}{10}\,
    \int\limits_{\ren} \hat{\psi}_2^*\,  y \cdot\nabla  \hat{\psi}_1,\quad 
   %%}
  %% \\ &
    %%\tex{
   E_2:=
    \frac{\a_2}{10}\, \int\limits_{\ren} \hat{\psi}_2^*\,  y
    \cdot (\nabla  \hat{\psi}_1-\nabla \hat{\psi}_2) -
    \hat{\psi}_3^*\,  y
    \cdot (\nabla \hat{\psi}_1 - \nabla\hat{\psi}_3).
   }
\end{align*}}

Therefore,  using the conic classification to solve \eqref{br64},
we have the number of solutions through the intersection of
two conics. Then, depending on the type of conic, we shall always
obtain one to four possible solutions for our system. Hence,
somehow, the number of solutions depends on the coefficients we
have for the system and, at the same time, on the eigenfunctions
that generate the subspace $\ker\big({\bf B}+\frac{k}{10}\big)$.

Thus, we have the following conditions, which will provide us with
the conic section of each equation of the system \eqref{br64}:
\begin{enumerate}
\item[(i)] If $B_i^2 - 4A_iE_i < 0$, the equation represents an {\em ellipse},
unless the conic is degenerate, for example $c_2^2 + c_3^2 + k = 0$
for some positive constant k. So, if $A_i=B_j$ and $E_i=0$,
the equation represents a {\em circle};
\item[(ii)] If $B_i^2 - 4A_iE_i = 0$, the equation represents a {\em parabola};
\item[(iii)] If $B_i^2 - 4A_iE_i > 0$, the equation represents a {\em hyperbola}.
If we also have $A_i + E_i = 0$ the equation represents a
hyperbola (a rectangular hyperbola).
\end{enumerate}
Consequently, the zeros of the system \eqref{br64} and, hence, of
the system \eqref{br61}, adding the ``normalizing" constraint
\eqref{br15}, are ascertained by the intersection of those two
conics in \eqref{br64} providing us with the number of possible
$n$-branches between one and {four}. Note that in case  those
conics are two circles we only have two intersection points at
most. Moreover, due to the dimension of the eigenspaces it looks
like in this  case that we have four possible intersection points
two of them will coincide. However, the justification for this is
far from clear.

Moreover, as was done for the
previous case when $|\b|=1$, we need to control the oscillations
of the perturbation functions in order to maintain the number of
solutions. Therefore, imposing that
\begin{equation*}
    \left\| \o_i (c_2,c_3) \right\|_{L^\infty} \leq \mf{F}_i(c_2^*,c_3^*),
    \quad \hbox{with}\quad i=1,2,
\end{equation*}
we ascertain that the number of solutions must be between one and four.
This again gives us an idea of the difficulty of more general
multiplicity results.

%%%%%%%%%%%%%%%%%%%%%%%%%%%%%%%%%%%%%%%%%%%%%%%%%
%%%%%%%%%%%%%%%%%%%%%%%%%%%%%%%%%%%%%%%%%%%%%%%%%%%%%%%%%
\section{Global extensions of bifurcation branches: numerical approach}
\label{S5}

Here we present numerical evidence for the nonlinear eigenfunctions whose eigenvalues are known explicitly. Namely the first eigenvalue-eigenfunction pair $\{\alpha_0(n),f_0\}$ and those in the $n=0$ case  $\{\alpha_k(0),f_k\}$. In these cases the eigenvalues are given explicitly by \eqref{alb1} and \eqref{bf4} respectively.

The first eigenvalue-eigenfunction pair $\{\alpha_0(n),f_0(|y|)\}$ satisfy \eqref{self1}, which may be integrated to 
\begin{equation}
\label{f0}
 \tex{
     |f_0|^{n} \frac{d}{d|y|} [ \D_y^4 f_0 ]  +  \a_0 |y| f_0 =0,\quad \hbox{where}\quad \D_y = \frac{d^2}{d|y|^2} + \frac{(N-1)}{|y|}\frac{d}{d|y|}, 
    }
\end{equation}
is the appropriate Laplacian radial operator. Use has been of the zero-flux and zero height conditions \eqref{i3} in self-similar form, which are imposed on the interface $|y|=y_0$ (i.e. $|x|=y_0 t^{\beta_0}$ with $\beta_0$ as given in \eqref{alb1}). Consequently, we add to \eqref{f0} the boundary conditions
 \begin{eqnarray}
 \mbox{at $y=0$:} \hspace{1cm} &&  f_0=1, \hskip 0.25cm \frac{d^{(i)} f_0}{d|y|^{(i)}}=0 \;\;\;\mbox{for $i=1,3,5,7$}, \label{f0bc1} \\
 \mbox{at $|y|=y_0$:} \hspace{1cm} &&  f_0=\frac{d^{(i)} f_0}{d|y|^{(i)}} =0 \;\;\;\mbox{for $i=1,2,3,4$}. \label{f0bc2} 
\end{eqnarray}
Since the $\alpha_0$ are known, this gives a tenth-order system when $n>0$ to determmine $f_0$ and the finite free boundary $y_0$. When $n=0$, then $y_0=\infty$. Figure \ref{fig1} shows illustrative $f_0$ profiles for selected $n$ values in one-dimension (N=1). The system was solved as an IVP in Matlab (shooting from $y=0$), using 
the ODE solver ode15s with error tolerances of AbsTol=RelTol=$10^{-10}$ and the regularisation $|f|^n= (f^2+\delta^2)^{n/2}$ with $\delta=10^{-10}$.

In the $n=0$ case, other eigenvalue-eigenfunction pairs $\{\alpha_k(0),f_k\}$ for $k\geq 1$ satisfy
\begin{equation}
\label{fkn0}
 \tex{
     \frac{d}{d|y|} \left[ \frac{d \D_y f_k}{d|y|} \right] + \frac{(N-1)}{|y|} \frac{d \D_y f_k}{d|y|} +  \frac{1}{10} |y| \frac{d f_k}{d|y|}   + \a_k(0) f_k =0,
    }
\end{equation}
with 
\begin{eqnarray}
 \mbox{at $y=0$:} \hspace{0.5cm} &&  \left\{ \begin{array}{ll} 
                             f_k=1, \hskip 0.25cm \frac{d^{(i)} f_k}{d|y|^{(i)}}=0 \;\;\;\mbox{for $i=1,3,5,7,9$}, \hspace{0.5cm} & \mbox{if $k$ is even}\label{fkn0bc1a} \\
                             \frac{d f_k}{d|y|}=1, \hskip 0.25cm f_k=\frac{d^{(i)} f_k}{d|y|^{(i)}}=0 \;\;\;\mbox{for $i=2,4,6,8$}, \hspace{0.5cm} &\mbox{if $k$ is odd} \label{fkn0bc1b} \\
                             \end{array} \right. 
\end{eqnarray}
and as $|y| \to \infty$:  $f_k \to 0$. 
%%\begin{eqnarray}
%%\hskip -9cm \mbox{as $|y| \to \infty$:} \hspace{0.5cm} &&  f_k \to 0. \label{fkn0bc2}
%%\end{eqnarray}
In regards to this last condition,
%\eqref{fkn0bc2}, 
we may determine from \eqref{fkn0} the actual asymptotic behaviour
\begin{equation}
   f_k \sim A |y|^{-\frac{4N}{9}} \exp \left( - \frac{9}{10}\alpha_k(0)^{\frac{1}{9}} \omega |y|^{\frac{10}{9}} \right), 
\end{equation}
for arbitrary constant A and $\omega$ may be a ninth root of unity $\omega^9=1$ with positive real part. This gives a five-dimensional stable bundle of asymptotic behaviours with
\[
     \omega=\exp\left( \pm\frac{2m\pi i}{9}\right), \hspace{1cm} m=0,1,2, 
\]
where the roots for $m=2$ have the smallest postive real parts and thus control the behaviour for large $|y|$. 
Figure \ref{fig2} show the eigenfunction profiles for the first four cases $k=0,1,2,3$, where the $k=0$ profile has been added and the same shooting numerical 
procedure used (appropriately adapted for this 10th-order system). The eigenfunctions have been arbitrarily normalised by $f_k(0)=1$ for $k$ even and $f_k'(0)=1$ for $k$ odd. 

The eigenvalue-eigenfunction pairs where the eigenvalues are not explicitly known, but have to be solved for, requires the solution of a 12th-order system. This will be discussed in \cite{AEG2}. 

%Deep numerical investigation of such difficult nonlinear eigenvalue problems is an inevitable step in modern nonlinear PDE theory, 
%and no specific references are needed to justify this obvious fact. It is clear that, for a long time (and, concerning some most delicate 
%results on the multiplicity of solutions and their spatial structures, probably never),
%any rigorous results will not be available.

%\com{JDE: do we need any preamble like that? Or just put numerics? do what you wish, just carefully explain the numerics and how it works. We will 
%add some math. around,  if necessary. This section should be rather short, but feel abs. free. If too much is obtained, we can move some parts 
%to the Paper-II, no problem.}

\begin{figure}[htp]

%\vskip -2.2cm
\hspace*{-1cm}
\includegraphics[scale=0.2]{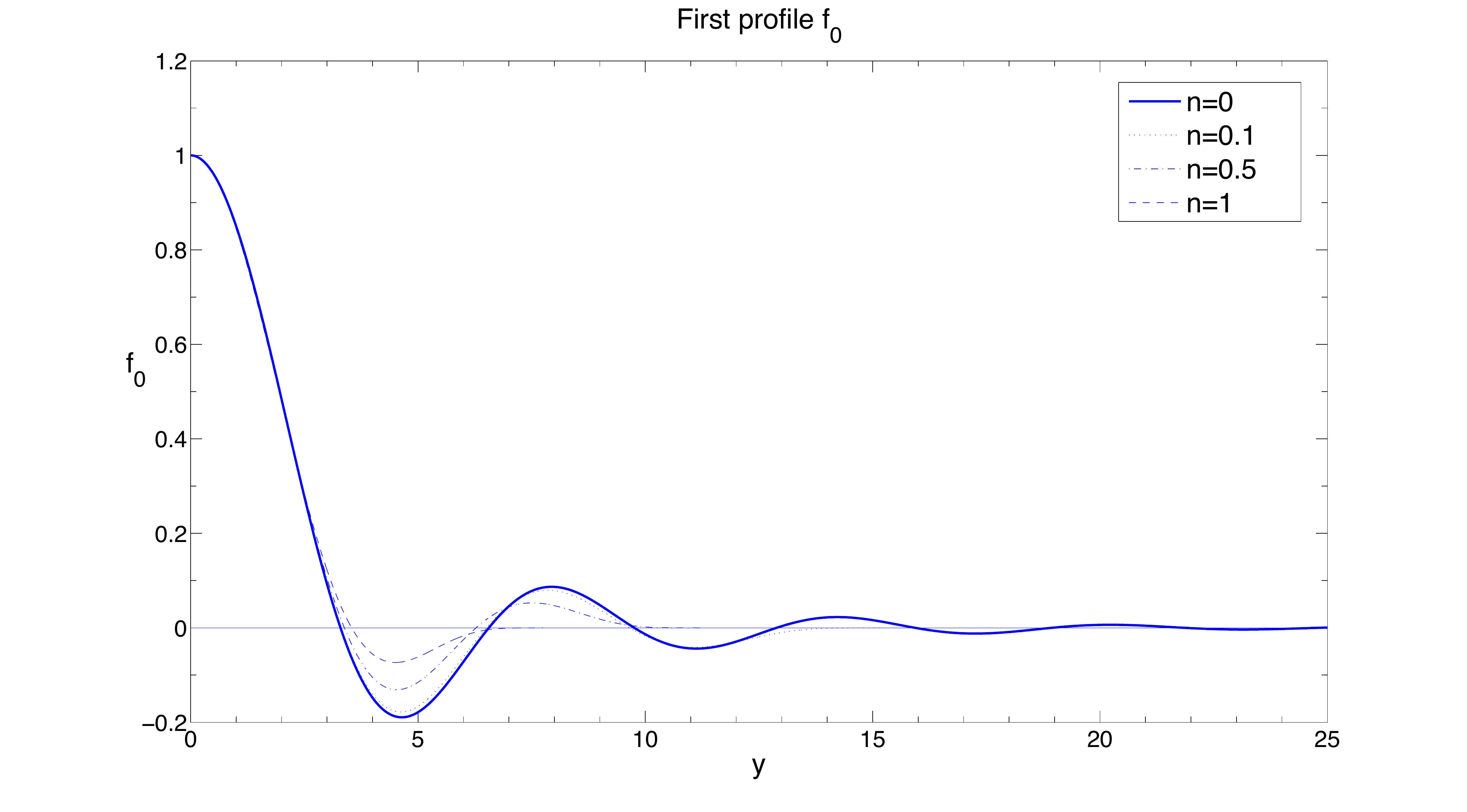}

\vskip -0.5cm \caption{ \small Profiles of the first eigenfunction $f_0$ for selected $n$. Obtained by numerical solution of \eqref{f0}--\eqref{f0bc2} in one-dimension $N=1$.  }
 \label{fig1}
%\end{center}
\end{figure}

\begin{figure}[htp]

%\vskip -2.2cm
\hspace*{-1cm}
\includegraphics[scale=0.2]{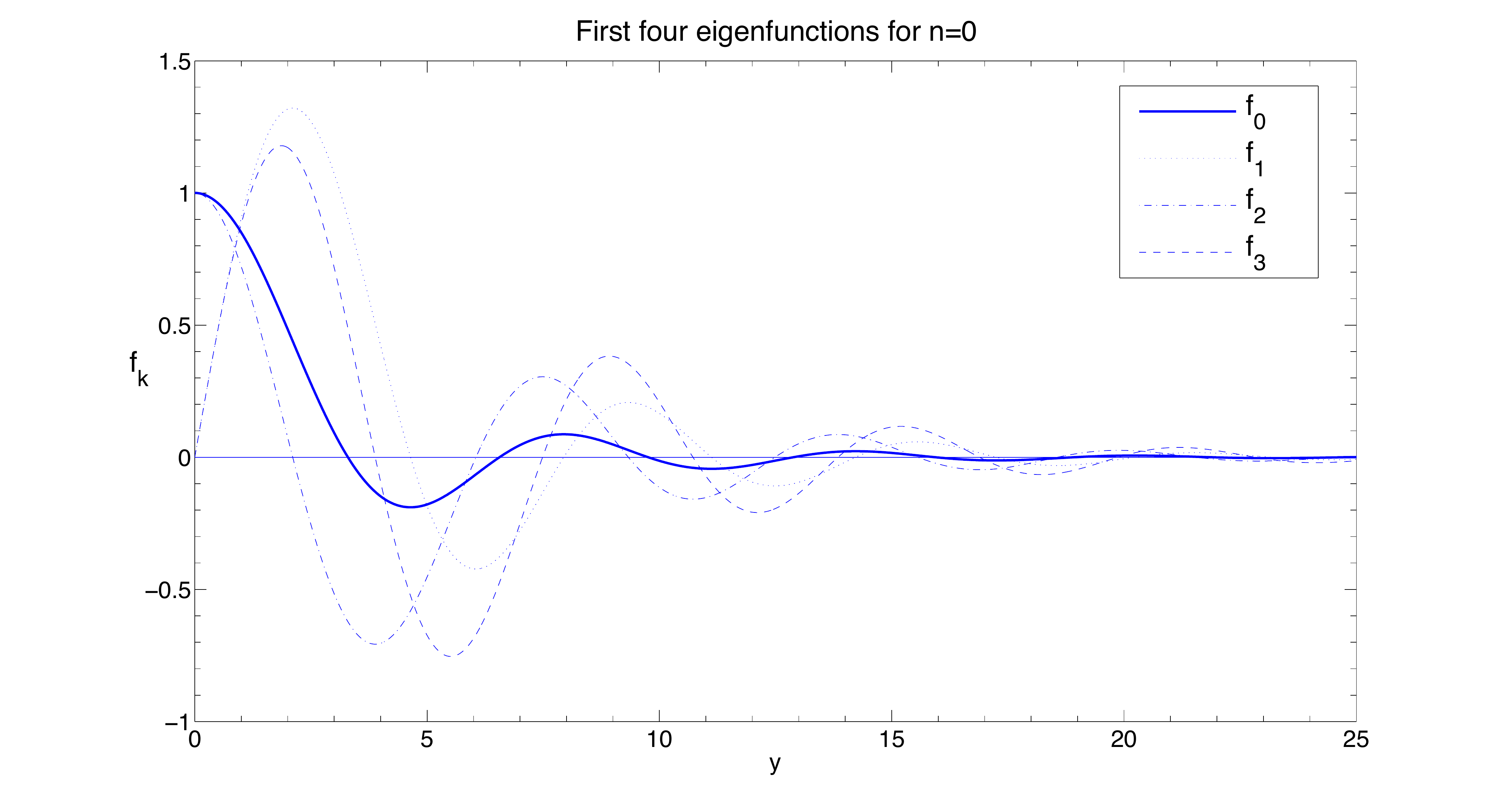}
%\includegraphics[width=7.3cm]{}
%\vskip -3.2cm \hspace{-1.5cm}
%\includegraphics[width=7.3cm]{}
%\includegraphics[width=7.3cm]{}

\vskip -0.5cm \caption{ \small Profiles of the first four eigenfunctions $f_k$, $k=0,1,2,3,$ in the case $n=0$. Obtained by numerical solution of \eqref{fkn0}--\eqref{fkn0bc2} in one-dimension $N=1$. }
 \label{fig2}
%\end{center}
\end{figure}

%%%%%%%%%%%%%%%%%%%%%%%%%%%%%%%%%%%%%%%%%%%%%%%%%%%%%%%%%%%%%%%%

%%%%%%%%%%%%%%%%%%%%%%%%%%%%%%%%%%%%%%%%%%%%%%%%
\begin{appendix}
\setcounter{section}{1} \setcounter{equation}{0}
\setcounter{subsection}{0}
\begin{small}

\section*{\bf Appendix A: Unstable TFE-10 model with an extra backward diffusion
term}
 \label{SState}

%%%%%%%%%%%%%%%%%%%%%%%%%%%%%%%%%%%%%%%%%%%%%%%%
%%%\section{Unstable TFE-10 model with an extra diffusion term}
 \label{Se1}

\subsection{Main model and problem setting}

\noindent
    Hereafter, we study the global-in-time behaviour of  solutions of
the tenth-order quasilinear evolution equation of parabolic type,
called the {unstable} TFE-10 \ef{e1},
%%\begin{equation}
%%\label{e1}
%%    u_{t} = -\nabla \cdot(|u|^{n} \n \D^4 u)-\D(|u|^{p-1}u)
%% \quad \mbox{in} \quad \ren \times \re_+
%%    \,,
%%\end{equation}
with the homogeneous  diffusion term of backward in time porous
medium type,
%% where $\n={\rm grad}_x$, $\D$ the usual Laplacian operator
 where $n>0$ and $p>n+1$ are given parameters. Equation \eqref{e1} is
also (as \eqref{i1}) written for solutions of changing sign, which
can occur in the CP and also in some FBPs.

%%%%%%%%%%%%%%%%%%%%%%%%%%%%%%%%%%%%%%%%%
%%\setcounter{equation}{0}
%%\section{Problem setting and self-similar solutions}
%%\label{Se2}

%%%%%%%%%%%%%%%%%%%%%%%%%%%%%%%%%%%

%%%%%%%%%%%%%%%%%%%%%%%%%%%%%%%%%%%
%%\subsection{The FBP and CP}

For  both the FBP and the CP, the solutions are assumed to satisfy
standard free-boundary conditions or boundary conditions at
infinity \eqref{i3} at the singularity surface (interface)
$\Gamma_0[u]$ given in \ef{gamma1}.
 For sufficiently smooth
interfaces,
%%%, for the FBP
the condition on the flux now reads
\begin{equation}
\label{flux}
    \lim_{\hbox{dist}(x,\Gamma_0[u])\downarrow 0}
    -{\bf n} \cdot  (|u|^{n} \nabla \Delta^4 u+ \nabla|u|^{p-1}u)=0.
\end{equation}
 Then,
differentiating  the mass $M(t)$ in \ef{Mass1} with respect to $t$
and applying the divergence theorem (under natural regularity
assumptions on solutions and free boundary), we get
\begin{equation*}
 \tex{
  J(t):=  \frac{{\mathrm d}M}{{\mathrm d}t}= -
  \int\limits_{\Gamma_0\cap\{t\}}{\bf n} \cdot 
     (|u|^{n} \nabla \Delta^4 u + \nabla |u|^{p-1}u)\, .
     }
\end{equation*}
The mass is conserved if
   $ J(t) \equiv 0$, which is assured by the flux condition
   \eqref{flux}.
The problem is completed with bounded, smooth, integrable,
compactly supported initial data denoted by \eqref{i4}.

In the CP for \eqref{e1} in $\ren \times \re_+$, one needs to pose
bounded compactly supported initial data \eqref{i4} prescribed in
$\ren$. Then,  under the same zero flux condition at finite
interfaces (to be established separately), the mass is preserved.

%%%%%%%%%%%%%%%%%%%%%%%%%%%%%%%%%%%%%%%%%%%%%%%%%%%%%
\subsection{Global similarity solutions}

 We now specify the self-similar
solutions of the equation \eqref{e1}, which are admitted due to
its natural scaling-invariant  nature. In the case of the mass
being conserved, we have global in time source-type solutions.
Using the following scaling in \eqref{e1}
%%\begin{equation*}
    $x:= \mu \bar x$,
    %\quad 
    $t:= \l \bar t$,
    %\quad 
    $u:= \nu \bar u$, 
    %\quad
%\end{equation*}
%%with
we obtain invariance provided 
%%\[
$\mu=\l^{\beta}$, 
%%\hskip 1cm  
$\nu = \l^{-\alpha}$,
%%\]
where 
\begin{equation}
\label{sfe3ab}
 \tex{
\a:=\frac{4}{5p-(n+5)} \quad \hbox{and}\quad \b:= \frac{p-(n+1)}{2[5p-(n+5)]}. 
    }
\end{equation}
This suggests considering similarity solutions of the form 
\begin{equation}
\label{sfe3}
 \tex{
    u(x,t):= t^{-\a} f(y), \quad \hbox{with}\quad
    y:=\frac{x}{t^\b}.\quad  
    }
\end{equation}
Substituting into \eqref{i1} and rearranging terms, we find that the function
$f$ solves a quasilinear elliptic equation of the form
\begin{equation}
\label{sfe4}
    \nabla \cdot \left[ |f|^{n} \n \D^4 f -\nabla (|f|^{p-1} f)\right] +\b \,y \cdot \nabla f+\a f=0\,.
\end{equation}
The parameters $\a$ and $\b$ (as given in \eqref{sfe3ab}) are linked by the following
expressions $$\tex{10\b-n\a=1,\quad 2\b-\a(p-1)=1.}$$

Finally, due to the above relations between $\a$ and $\b$, we find a {\em
nonlinear eigenvalue problem} of the form
\begin{equation}
\label{sfe5}
 \fbox{$
 \tex{
     \nabla \cdot \left[ |f|^{n} \n \D^4 f -\nabla (|f|^{p-1} f)\right] +\frac{1+\a n}{10}\, y \cdot \nabla f +\a
    f=0,
    \quad f \in C_0(\ren)\, ,
    }
    $}
\end{equation}
 where we add to the equation \ef{sfe4} a natural assumption that
 $f$ must be compactly supported (and, of course, sufficiently
 smooth at the interface, which is an accompanying question to be
 discussed as well).

Thus, for such degenerate elliptic equations,
  the functional setting in \ef{sfe5} assumes that we are
 looking for  (weak) {\em compactly supported} solutions $f(y)$ as
 certain ``nonlinear eigenfunctions" that hopefully occur for special values of nonlinear eigenvalues
  $\{\a_\g\}_{|\g| \ge 0}$.
 Similar to the previous problem, we intend to justify (formally,
 at least) that \ef{main1} holds for the problem \ef{sfe5}.
%%   Our goal is to justify
%%  that, labelling  the eigenfunctions via a multiindex $\g$,
%% \be
%% \label{sfe51}
%% \fbox{$
%% \mbox{
%% (\ref{sfe5}) possesses  a countable set of
%% eigenfunction/value pairs $\{f_\g,\, \a_\g\}_{|\g|=k \ge 0}$.
%%  }
%%  $}
%%  \ee
Moreover, again for this particular situation, in the linear case
$n=0$,
 the condition $f \in C_0(\ren)$, is  replaced by the requirement that the
 eigenfunctions $\psi_\b(y)$ exhibit typical exponential decay at
 infinity by using the weighted space \ef{WW11}.

Next, using the mass evolution \ef{Mass11}, in the case $\int f
\not = 0$,
  the exponents are calculated giving the first explicit nonlinear eigenvalue:
\begin{equation}
\label{albe1}
 \tex{
 -\a + \b N=0 \LongA
   p_0(n)=n+1+\frac{8}{N},\quad  \a_0(n)=\frac{N}{10+Nn} \quad \mbox{and}  \quad \b_0(n)=\frac{1}{10+Nn}.
    }
\end{equation}

So far, the analysis looks rather similar to the one performed
previously for the 10$th$--order equation without the extra
diffusion term \eqref{self1}. However, the results seem to be
quite different. The main difference is that, for  \eqref{self1}
(rescaled version of \eqref{i1}), we ascertained the
branching--asymptotic analysis from the solutions or
eigenfunctions of the rescaled poly-harmonic equation \eqref{s4}.
 For     \eqref{sfe5}, the solutions will emanate from the solutions
of a nonlinear perturbation of the equation \eqref{s4}, basically
due to the extra diffusion term.

It was obtained in \cite{EGK2} that, for the fourth--order
unstable TFE4,
\begin{equation}
\label{foe4}
    u_{t} = -\nabla \cdot(|u|^{n} \nabla\Delta u)-\D(|u|^{p-1}u)
 \quad \mbox{in} \quad \ren \times \re_+
    \,,
\end{equation}
there are continuous families of solutions of global similarity
solutions when the exponent $p$ is the critical exponent $
 \tex{
 p=p_0=
n+1+\frac{4}{N}.
 }
 $
  Moreover, the authors also showed that in the
particular case when $p\neq p_0$ the families of similarity
solutions become countable.

Let us briefly comment on that. Namely, in 1D,
 the main reason in the critical case $p=p_0$ to admit wider (a
 continuum) family of solutions is that the corresponding rescaled
 ODE admits integration once and reduces to a {\em third-order}
 ODE, which makes a shooting procedure underdetermined: two parameters to satisfy a single
 symmetry conditions at the origin.
  For $p \neq p_0$, the ODE is truly {\em fourth-order}, and the shooting
 is well-posed: two parameters and two symmetry conditions.

\ssk

 A similar situation occur the above unstable TFE-10: for
$p=p_0$ there exists symmetry reduction and the ODE in 1D becomes
of ninth order.
 This analysis could be extended  to
our 10$th$--order equation Therefore, \eqref{e1} admits continuous
families of global similarity solutions if $p=p_0$ given in
\ef{albe1}
%% the critical exponent, such that $$
%% p_0= n+1+\frac{8}{N},$$
and, also, for
$p\neq p_0$ we will have a countable family of solutions for the
unstable TFE--10 \eqref{e1}.

For equations in $\ren$, a similar result holds true in the radial
setting, where we deal with ODEs again. Non-radial patterns are
entirely unknown and, honestly, we do not have any clue how and by
what tools these can be detected (numerics are expected also to be
extremely difficult).

Therefore (in the in ODE setting\footnote{This requirement will shortly be relaxed; see below.}), performing a similar branching analysis, as
the one done in the previous section for the TFE--10 \eqref{i1},
we obtain that \eqref{e1} possesses  a countable set of
eigenfunction/value pairs $\{f_k,\, \a_k\}_{k \ge 0}$
\eqref{main1} such that the solutions of the equation \eqref{sfe5}
emanate from the solutions of the rescaled version Cahn--Hilliard
equation type
 \begin{equation}
\label{chtp}
    u_{t} = \Delta^5 u-\D(|u|^{p-1}u)
 \quad \mbox{in} \quad \ren \times \re_+,
    \,,
\end{equation}
at $n=0$. In other words, the solutions of the equation
\begin{equation}
\label{recah}
 %%\fbox{$
 \tex{
     \Delta^5 f -\Delta (|f|^{p-1} f)+\frac{1}{10}\, y \cdot \nabla f +\a
    f=0,
    \quad f \in H^{10}_\rho(\ren)\, ,
    }
  %%  $}
\end{equation}
for certain values of the parameter $\a$, which will provide us
with that countable family of solutions emanating form the
solutions of \eqref{recah} at $n=0$. One can easily see that
\eqref{recah} is a nonlinear perturbation of the rescaled equation
\eqref{s4}.
 Moreover, to detect a deeper connection with linear
 eigenfunctions, a further homotopy deformation analysis
 should be performed by passing $p \to 1^+$ leading to the linear
 eigenvalue problem
\begin{equation}
\label{recahN}
 %%\fbox{$
 \tex{
     \Delta^5 f -\Delta f +\frac{1}{10}\, y \cdot \nabla f +\a
    f=0,
    \quad f \in H^{10}_\rho(\ren)\, ,
    }
  %%  $}
\end{equation}
which admits a clear study similar to \cite{EGKP}. It is important
that we can describe the whole complete family of eigenfunctions
of \ef{recahN} including all the non-radial ones.

Thus, it turns out that the solutions of the equation \eqref{e1}
can emanate from a nonlinear perturbed version of the
eigenfunctions for the equation \eqref{s4} via two-parametric
homotopy deformation to a linear eigenvalue problem. This, at
least, very formally explains the origin of countablity of
nonlinear eigenfunctions family of those TFEs-10.

\end{small}
\end{appendix}

\end{document}